\documentclass[10pt]{amsart}

\usepackage[utf8]{inputenc}
\usepackage[english]{babel}
\usepackage{amsmath}
\usepackage{float}
\usepackage{amssymb}
\PassOptionsToPackage{dvipsnames}{xcolor}
\usepackage{tikz}
\usepackage{esint}
\usepackage{amsfonts}
\usepackage{multicol}

\newtheorem{innercustomthm}{Theorem}
\newenvironment{customthm}[1]
  {\renewcommand\theinnercustomthm{#1}\innercustomthm}
  {\endinnercustomthm}
  
 \newtheorem{innercustomcorollary}{Corollary}
\newenvironment{customcorollary}[1]
  {\renewcommand\theinnercustomcorollary{#1}\innercustomcorollary}
  {\endinnercustomcorollary}

  \newtheorem{innercustomlemma}{Lemma}
\newenvironment{customlemma}[1]
  {\renewcommand\theinnercustomlemma{#1}\innercustomlemma}
  {\endinnercustomlemma}

\newtheorem*{theorem*}{``Theorem''}

\newcommand{\ts}{\sigma^{\h}}
\newcommand{\s}{\sigma}

\newcommand{\tp}{\phi^{\h}}

\newcommand{\pj}{\partial_j}
\newcommand{\pk}{\partial_k}
\newcommand{\df}{\displaystyle\frac{1}{r}}

\newcommand{\di}{\nabla \cdot}
\newcommand{\hs}{\mathbb{H}^d_+}
\newcommand{\h}{\mathbb{H}}
\newcommand{\bhs}{\partial \mathbb{H}^d_+ }
\newcommand{\cd}{\nabla \tp_{d,M} - \nabla \phi_{d}}

\newcommand{\ir}[2]{ \left( \displaystyle\int_{B_{#1}^+}|{#2}|^2 \,dx\right)^{1/2}}
\newcommand{\irs}[2]{ \displaystyle\int_{B_{#1}^+}|{#2}|^2 \,dx}

\newcommand{\irt}[2]{ \displaystyle\int_{#1}|{#2}|^2 \,dx}
\newcommand{\irtr}[2]{ \left( \displaystyle\int_{#1}|{#2}|^2 \,dx\right)^{1/2} }

\newcommand{\ipts}[2]{ \displaystyle\int_{#1}{#2}\,dx}
\newcommand{\fir}[2]{ \left( \displaystyle\fint_{B_{#1}^+}|{#2}|^2 \,dx\right)^{1/2}}

\newcommand{\firs}[2]{  \displaystyle\fint_{B_{#1}^+}|{#2}|^2 \,dx}

\newcommand{\Rp}{R^{\prime}}
\begin{document}

\title[Random Elliptic Operators on the Half-Space]{Liouville Principles and a Large-Scale Regularity Theory for Random Elliptic Operators on the Half-Space}
\author{Julian Fischer}
\author{Claudia Raithel}
\begin{abstract}
We consider the large-scale regularity of solutions to second-order linear elliptic equations with random coefficient fields. In contrast to previous works on regularity theory for random elliptic operators, our interest is in the regularity at the boundary: We consider problems posed on the half-space with homogeneous Dirichlet boundary conditions and derive an associated $C^{1,\alpha}$-type large-scale regularity theory in the form of a corresponding decay estimate for the homogenization-adapted tilt-excess. This regularity theory entails an associated Liouville-type theorem. The results are based on the existence of homogenization correctors adapted to the half-space setting, which we construct -- by an entirely deterministic argument -- as a modification of the homogenization corrector on the whole space. This adaption procedure is carried out inductively on larger scales, crucially relying on the regularity theory already established on smaller scales.
\end{abstract}

\maketitle

\section{Introduction}

Classical counterexamples in the theory of the second-order linear elliptic equation
\begin{align}
\label{LinearEllipticEquation}
-\nabla \cdot (a\nabla u) = 0
\quad\text{on }\mathbb{R}^d
\end{align}
demonstrate that uniform ellipticity and boundedness of the coefficient field $a$ are not sufficient to ensure Lipschitz continuity of weak solutions: It is well-known that for any H\"older exponent $0<\alpha\leq 1$ there exists a uniformly elliptic coefficient field $a$ and an associated weak solution $u\in H^1_{loc}(\mathbb{R}^d)$ which fails to be H\"older continuous with exponent $\alpha$ (see e.\,g.\ the example of Meyers \cite[Example 3]{PiccininiSpagnolo}). For second-order linear elliptic systems, the celebrated counterexample of De~Giorgi (see e.\,g.\ \cite[Section 9.1.1]{GiaquintaMartinazzi}) shows that weak solutions may even locally fail to be bounded.

H\"older regularity of a function $u$ is equivalent to suitable approximability of $u$ by polynomials, a property that may be rephrased in terms of a decay estimate for the associated tilt-excess. An interesting consequence of quantitative regularity estimates on large scales are Liouville principles: Liouville principles provide a characterization of the dimension of the space of solutions on $\mathbb{R}^d$ which satisfy a given polynomial growth restriction. In fact, the abovementioned classical counterexamples to regularity are at the same time counterexamples to the corresponding Liouville principles.

These classical counterexamples to regularity share the feature of imposing a certain large-scale structure on the coefficient field $a$: For example, both in the counterexample of Meyers and in the counterexample of De~Giorgi the coefficient field $a$ has a purely radial structure. In a series of recent works \cite{BenjaminiCopinKozmaYadin,MarahrensOtto,ArmstrongSmart,ArmstrongMourrat,GloriaNeukammOtto,FischerOtto}, it has been established that the coefficient fields which constitute such counterexamples are in fact in a certain sense necessarily ``non-generic'': For \emph{random} coefficient fields -- more precisely, coefficient fields chosen according to a stationary and ergodic probability measure on the space of uniformly elliptic and bounded coefficient fields (see below for a definition) -- , almost surely a large-scale regularity theory in the form of a corresponding decay estimate for the tilt-excess on large scales and Liouville principles hold.

Motivated by these recent results, in the present work we consider the large-scale \emph{boundary regularity} of solutions to linear elliptic equations with random coefficient fields in the case of Dirichlet boundary conditions. More precisely, we consider solutions to the problem
\begin{subequations}
\label{firstequation}
\begin{align}
 -\di (a \nabla u) &= 0 \quad \textrm{ in } \hs,\\
  u &\equiv 0 \quad \textrm{ on } \bhs,
\end{align}
\end{subequations}
with $a$ being the restriction of a random coefficient field on the full space $\mathbb{R}^d$ to the half-space
\begin{align*}
\hs:=\{(x_1,\ldots,x_{d})\in \mathbb{R}^d: x_d > 0\}.
\end{align*}
The main results of the present work are a ``$C^{1,\alpha}$-type regularity theory on large scales'' for such solutions -- in the form of a $C^{1,\alpha}$-type excess-decay estimate for the tilt-excess on large scales -- and an associated Liouville-type theorem.

Randomness in the coefficient field does not exclude the possible occurrence of counterexamples to regularity on small scales. For this reason, one may only hope to establish an improved regularity theory for such random elliptic operators on large scales. A rigorous mathematical meaning to the notion of ``large-scale regularity'' may be given in terms of corresponding decay estimates for the tilt-excess: 
The classical notion of tilt-excess compares a solution of the elliptic equation $-\nabla \cdot(a\nabla u)=0$ to e.\,g.\ the space of affine polynomials $x\mapsto \xi\cdot x+c$ in the squared energy norm. For a function $u$, the tilt-excess on the ball $\{|x|<r\}$ is defined as
\begin{align*}
\operatorname{Exc}(r):=\inf_{\xi\in \mathbb{R}^d} \fint_{\{|x|<r\}} |\nabla u-\xi|^2 \,dx.
\end{align*}
Differentiability properties of the function $u$ are then encoded in decay properties of the tilt-excess in the radius $r$. For solutions to the Laplace equation $-\Delta u=0$ on $\mathbb{R}^d$, the tilt-excess displays decay in the radius $r$ according to
\begin{align*}
\operatorname{Exc}(r)\leq \left(\frac{r}{R}\right)^2 \operatorname{Exc}(R)
\end{align*}
for any pair of radii $0< r\leq R$. When valid for balls $\{|x-x_0|<r\}$ around any center $x_0\in \mathbb{R}^d$, this excess-decay property entails $C^{1,1}$-regularity of solutions. It is now intuitive to give a meaning to the notion of large-scale regularity of a function $u$ by asking for appropriate excess-decay on large scales, i.\,e.\ excess-decay for radii larger than a certain minimal radius $r^\ast$.

This classical definition of tilt-excess is, however, not the appropriate quantity in the framework of random coefficient fields $a$ on $\mathbb{R}^d$: In this case, one does not expect the fluctuations of $\nabla u$ around some constant value $\xi$ to be small. It is therefore necessary to suitably adapt the notion of tilt-excess to this setting. In \cite{GloriaNeukammOtto}, motivated by a similar ansatz of Avellaneda and Lin \cite{AvellanedaLinCPAM} in the context of periodic homogenization, Gloria, Neukamm, and Otto have introduced the \emph{homogenization-adapted} notion of tilt-excess
\begin{align*}
\operatorname{Exc}(r):=\inf_{\xi\in \mathbb{R}^d} \fint_{\{|x|<r\}} |\nabla u-\nabla (\xi\cdot x+\phi_\xi)|^2 \,dx.
\end{align*}
Here, $\phi_\xi$ is the so-called \emph{homogenization corrector} (see below for a definition). We would like to emphasize that the ``corrected affine polynomial'' $\xi\cdot x+\phi_\xi$ may be regarded as a perturbation of the original polynomial $\xi\cdot x$, while at the level of the gradient $\xi+\nabla \phi_\xi$ is typically \emph{not} a perturbation of $\xi$.

In the half-space setting, the Dirichlet boundary conditions on $\bhs$ introduce further restrictions on solutions $u$ of the problem \eqref{firstequation}: It will turn out that the solutions may be approximated in terms of a multiple of just the perturbed coordinate function $x_d+\tp_d$ (with $\tp_d$ denoting the homogenization corrector adapted to the Dirichlet boundary condition on $\bhs$), lifting the need to consider perturbations of general affine functions $\xi\cdot x$. The appropriate notion of tilt-excess in our setting of the equation \eqref{firstequation} is therefore given by the formula \eqref{OurTiltExcess} below. Omitting for the moment the precise assumptions on the random coefficient field, our main result with respect to the large-scale regularity of solutions to the problem \eqref{firstequation} may be phrased as follows:
\begin{theorem*}
Define the tilt-excess of a function $u$ on the half-ball
\begin{align*}
B_r^+:=\{|x|<r\}\cap \hs
\end{align*}
as
\begin{align}
\label{OurTiltExcess}
\textrm{Exc}^{\h}(r):=\inf_{b\in \mathbb{R}} \fint_{B_r^+} |\nabla u-b(e_d+\nabla\tp_d)|^2 \,dx.
\end{align}
Let $0<\alpha<1$ and let $a$ be a random coefficient field subject to our assumptions on the random coefficient field stated below. Then almost surely, the following holds:
\begin{itemize}
\item[i)] There exists a homogenization corrector $\phi_d^\h$ which solves the corrector equation
\begin{subequations}
\begin{align}
-\nabla \cdot (a\nabla \phi_d^\h) &=\nabla \cdot (ae_d) &&\text{in }\quad\hs,
\\
\phi_d^\h &=0 &&\text{on }\quad\bhs,
\end{align}
\end{subequations}
and satisfies the sublinear growth condition
\begin{align*}
\lim_{r\rightarrow \infty} \frac{1}{r} \bigg(\fint_{B_r^+} |\phi_d^\h|^2 \,dx\bigg)^{1/2} =0.
\end{align*}
\item[ii)] There exists a finite $r^\ast$ such that any weak solution to the problem \eqref{firstequation} satisfies the excess-decay estimate
\begin{align*}
\textrm{Exc}^{\h}(r) \lesssim \left(\frac{r}{R}\right)^{2\alpha}\textrm{Exc}^{\h}(R)
\end{align*}
for any pair of radii $R\geq r\geq r^\ast$.
\end{itemize}
\end{theorem*}
The classical (zeroth-order) Liouville theorem states that any solution $u$ to the Laplace equation $-\Delta u=0$ on $\mathbb{R}^d$ with sublinear growth at infinity must be constant. More generally, any harmonic function on $\mathbb{R}^d$ which satisfies a growth condition of the form $|u(x)|\lesssim 1+|x|^{k+\alpha}$ (with $k\in \mathbb{N}_0$, $0<\alpha<1$) is equal to a harmonic polynomial of order less or equal to $k$.
For the Laplacian on the half-space with homogeneous Dirichlet boundary conditions, half-space-adapted Liouville principles are available: The first-order Liouville principle states that any solution $u$ to the equation $-\Delta u=0$ on $\mathbb{H}^d_+$ with $u\equiv 0$ on $\partial \mathbb{H}^d_+$ and with subquadratic growth in the sense
\begin{align}
\label{Subquadratic}
|u(x)|\lesssim 1+|x|^{1+\alpha}
\end{align}
for some $\alpha<1$ is a multiple of the coordinate function
\begin{align*}
x\mapsto x_d.
\end{align*}
In the present work, for random coefficient fields we shall similarly characterize the subquadratically growing solutions $u$ to the equation $-\nabla \cdot(a \nabla u)=0$ with homogeneous Dirichlet boundary conditions on $\bhs$:
\begin{theorem*}
Let $a$ be a random coefficient field subject to our assumptions on the random coefficient field stated below. Then almost surely, the following assertion holds: Any weak solution $u$ to the problem \eqref{firstequation} which satisfies a growth condition of the form
\begin{align*}
\lim_{r\rightarrow \infty} \frac{1}{r^{1+\alpha}}\left(\fint_{B_r^+} |u|^2 \,dx \right)^{1/2} =0
\end{align*}
for some $\alpha\in (0,1)$ is a multiple of the ``perturbed coordinate function''
\begin{align*}
x\mapsto x_d+\tp_d(x).
\end{align*}
\end{theorem*}
This Liouville principle is a simple consequence of the excess-decay estimate on large scales.

For the results of our present work, by a \emph{random coefficient field} we shall understand a coefficient field chosen at random according to some probability measure $\langle\cdot\rangle$ (which is also called ``ensemble'' in this context) on the space of coefficient fields on $\mathbb{R}^d$. Our two main assumptions on the ensemble are the following:
\begin{itemize}
\item The assumption of stationarity (shift-invariance), stating that the measure $\langle \cdot \rangle$ is invariant under simultaneous spatial translation of all coefficient fields.
\item The assumption of ergodicity, which requires that any shift-invariant random variable must be $\langle\cdot\rangle$-almost surely constant, corresponding to a qualitative assumption on decorrelation on large scales. In the present work, we will need a slightly strengthened (slightly quantified) version of qualitative ergodicity, expressed in form of the growth estimate for the corrector \eqref{condition} below.
\end{itemize}
In addition, we shall assume that the probability measure $\langle \cdot \rangle$ is supported on uniformly elliptic and bounded coefficient fields: We require that there exists a constant $\lambda>0$ such that almost surely for almost every $x\in \mathbb{R}^d$ the estimates
\begin{subequations}
\label{Elliptic}
\begin{align}
|a(x)v|&\leq |v|,
\\
a(x)v \cdot v &\geq \lambda |v|^2 
\end{align}
\end{subequations}
hold for every $v\in \mathbb{R}^d$. Note that the choice $|a(x)v|\leq |v|$ for the upper bound is out of convenience and does not lead to a loss of generality, as given a general upper bound it may be enforced by rescaling.

\begin{figure}[H]
\begin{tikzpicture}[scale=0.3]
\draw[fill=Blue,color=Blue] (0,0) rectangle (10,10);
\clip (0,0) rectangle (10,10);
\foreach \x in {(0.86,8.91),(2.88,9.83),(0.34,0.45),(5.53,7.27),(1.16,7.84),(1.86,2.03),(4.09,0.80),(2.71,9.88),(7.60,2.98),(3.22,6.60),(4.19,2.04),(9.38,1.23),(2.52,4.79),(9.86,8.18),(1.59,7.22),(5.63,5.96),(1.96,4.86),(9.27,4.12),(0.72,7.81),(4.33,9.54),(1.33,1.73),(6.02,8.74),(5.70,6.52),(3.45,0.23),(8.39,3.82),(4.18,0.95),(0.57,7.01),(9.56,8.30),(0.21,5.96),(5.37,2.03),(3.43,3.12),(8.45,9.17),(5.70,1.31),(4.90,1.67),(5.80,4.88),(7.33,4.00),(9.57,1.10),(5.60,2.61),(1.53,7.89),(3.18,3.82),(7.39,8.41),(9.93,9.56),(9.13,6.67),(2.38,4.77),(3.26,7.11),(2.71,7.92),(1.66,6.58),(5.70,6.94),(0.36,5.59),(1.22,5.97),(9.13,2.12),(3.89,7.63),(7.05,6.38),(2.34,0.49),(3.99,6.62),(0.41,5.40),(6.45,0.98),(1.72,2.65),(1.63,4.15),(1.45,5.21),(5.60,5.92),(9.42,9.68),(7.25,4.06),(5.20,5.43),(8.29,7.89),(5.06,0.86),(5.03,6.95),(6.91,7.39),(7.60,2.05),(4.61,9.36),(2.60,3.32),(7.99,5.65),(4.35,8.96),(0.92,8.48),(6.67,4.37),(1.49,5.95),(6.73,8.80),(5.45,5.09),(7.07,6.34),(3.78,6.80),(6.20,4.44),(8.51,5.97),(1.98,8.27),(5.29,8.83),(3.58,0.80),(8.39,7.65),(5.09,4.50),(6.06,5.45),(5.82,3.57),(6.90,0.07),(1.91,6.37),(7.13,8.55),(9.09,9.07),(5.11,2.78),(0.88,0.92),(4.81,6.71),(9.93,9.50),(0.93,3.87)
}
\draw[fill=Maroon,color=Maroon] \x circle (0.5);
\end{tikzpicture}
~
\includegraphics[scale=0.285]{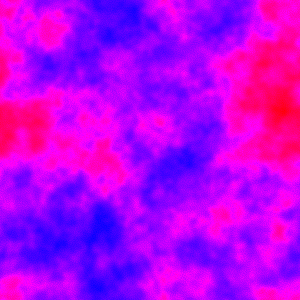}
\caption{Two examples of random coefficient fields.\label{FigureRandom}}
\end{figure}
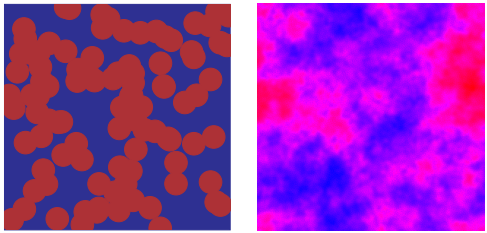

To name a few examples, our results in the present work are in particular applicable to the following cases of ensembles of random coefficient fields:
\begin{itemize}
\item Ensembles for which $a(x)$ is either equal to a positive definite matrix $a_1$ or equal to another positive definite matrix $a_2$, depending on whether $x$ is contained in a random set of balls of a given fixed radius, the centers of the balls being chosen according to a Poisson point process (see the left picture in Figure \ref{FigureRandom}).
\item Stationary ensembles with finite range of dependence (i.\,e.\ ensembles for which $a|_U$ and $a|_V$ are stochastically independent for any two sets $U,V\subset\mathbb{R}^d$ with $\operatorname{dist}(U,V)\geq c$) subject to uniform ellipticity and boundedness conditions. Note that the previous case is a particular case of this.
\item Coefficient fields of the form $\xi(\tilde a(x))$, where $\tilde a$ denotes a matrix-valued stationary Gaussian random field subject to the decorrelation estimate
\begin{align*}
|\operatorname{Cov}(\tilde a(x),\tilde a(y))|\leq \frac{C}{|x-y|^\beta}
\end{align*}
for some $\beta \in (0,d)$ and where $\xi:\mathbb{R}^{d\times d}\rightarrow \mathbb{R}^{d\times d}$ is a Lipschitz map taking values in a bounded uniformly elliptic subset of the matrices of dimension $d\times d$ (see the right picture in Figure \ref{FigureRandom}).
\end{itemize}
That our results apply to the second example -- i.\,e.\ that condition \eqref{condition} below is satisfied almost surely for such ensembles -- follows e.\,g.\ from the estimates in \cite{GloriaOttoNew}. That our results are applicable to the third example is shown in \cite{CorrectorEstimatesSlowDecorrelation}.

Generally speaking, the improvement in the regularity of solutions to elliptic equations with random coefficient fields on large scales may be viewed as a homogenization effect: Classical results in qualitative stochastic homogenization state that on large scales the behavior of the second-order linear elliptic equation with a random coefficient field is close to the behavior of a constant-coefficient equation \cite{PapanicolaouVaradhan}. In fact, Avellaneda and Lin had established Liouville-type theorems \cite{AvellanedaLin} and regularity results \cite{AvellanedaLinCPAM} for \emph{periodic} coefficient fields -- i.\,e.\ in the context of \emph{periodic homogenization} -- long before the first works on random coefficient fields.

To the best of our knowledge, the first result on improved large-scale regularity properties of random elliptic operators has been derived by Benjamini, Duminil-Copin, Kozma, and Yadin \cite{BenjaminiCopinKozmaYadin} in the form of a zeroth-order Liouville theorem in the context of random walks in random environments. Their result holds under the assumptions of stationarity and qualitative ergodicity and includes the case of percolation, i.\,e.\ also suitable coefficient fields which are not uniformly elliptic.
In the work of Marahrens and Otto \cite{MarahrensOtto}, a large-scale $C^{0,\alpha}$-type regularity theory for any $\alpha<1$ was developed, assuming a quantification of ergodicity in the form of a logarithmic Sobolev inequality. In the work of Armstrong and Smart \cite{ArmstrongSmart}, a large-scale $C^{0,1}$-type regularity theory has been established under the assumption of finite range of dependence. Motivated by \cite{ArmstrongMourrat}, Gloria, Neukamm, and Otto \cite{GloriaNeukammOtto} derived a large-scale $C^{1,\alpha}$-type regularity theory in the form of a corresponding excess-decay estimate and a Liouville principle; their result is applicable in the case of just stationarity and qualitative ergodicity of the ensemble. Finally, the picture of large-scale regularity was mostly completed by Otto and the first author \cite{FischerOtto}, who developed a $C^{k,\alpha}$-type large-scale regularity theory and associated Liouville principles, assuming only a slight quantification of ergodicity. Later, another proof for such a large-scale $C^{k,\alpha}$-type regularity theory was given by Armstrong, Kuusi, and Mourrat \cite{ArmstrongKuusiMourratComplete}; while the results of \cite{ArmstrongKuusiMourratComplete} are stated under the assumption of finite range of dependence of the ensemble, as also mentioned in \cite{ArmstrongKuusiMourratComplete} it is apparent from their proof and \cite{ArmstrongMourrat} that their arguments also apply to settings with weak quantitative decorrelation. Recently, a large-scale $C^{k,\alpha}$ regularity theory for random elliptic operators on Bernoulli percolation clusters has been developed by Armstrong and Dario \cite{ArmstrongDario}.

The \emph{periodic} (and almost periodic) homogenization of boundary value problems for the elliptic equation $-\nabla \cdot(a\nabla u)=f$ has a long history. Avellaneda and Lin \cite{AvellanedaLinCPAM} have derived a $C^{0,1}$ regularity theory up to the boundary in the context of periodic homogenization with Dirichlet boundary conditions. In their work, they also adapt the homogenization correctors to the Dirichlet boundary conditions, however only locally and in a different way on every scale. For Neumann boundary conditions, the corresponding regularity result has been obtained by Kenig, Lin, and Shen \cite{KenigLinShenNeumann}. In the almost periodic case, a $C^{0,1}$ regularity theory for Dirichlet and Neumann boundary conditions has been established by Armstrong and Shen \cite{ArmstrongShen}; though not mentioned in the paper, their arguments -- reminiscent of the ones in \cite{ArmstrongMourrat} -- likely also apply to the setting of stochastic homogenization.

Rates of convergence for the periodic homogenization of elliptic equations on bounded domains have been established by Avellaneda and Lin \cite{AvellanedaLinCPAM} and Kenig, Lin, and Shen \cite{KenigLinShen} in the Dirichlet and Neumann case, respectively. Higher-order approximations for periodic homogenization problems on bounded domains with Dirichlet boundary conditions have been obtained via boundary layer correctors by Allaire and Amar \cite{AllaireAmar} in the case of polygonal domains with rational slopes; G\'erard-Varet and Masmoudi \cite{GerardVaretMasmoudi2} have treated the case of polygonal domains with diophantine normals.

Basically, the introduction of boundary layer correctors leads to a homogenization problem with oscillating boundary data.
%by Kenig, Lin, and Shen \cite{KenigLinShen} in the case of boundary data with bounded Sobolev norm, even for just Lipschitz domains.
In the case of oscillating Dirichlet boundary data on general (even smooth) domains, the convergence properties may be drastically worse compared to the case of smooth boundary data: A recent result by Aleksanyan \cite{Aleksanyan} shows that the convergence may be arbitrarily slow. For uniformly convex domains, rates of convergence may be obtained as shown by G\'erard-Varet and Masmoudi \cite{GerardVaretMasmoudi}.
In the recent work by Armstrong, Kuusi, Mourrat, and Prange \cite{ArmstrongKuusiMourratPrange}, improved convergence rates have been derived; for $d\geq 4$, their rates reach the optimal exponent from the model case of constant coefficients (a case that was treated by Aleksanyan, Shahgholian and Sj\"olin \cite{AleksanyanShahgholianSjolin}). The subsequent improvement of \cite{ArmstrongKuusiMourratPrange} by Shen and Zhuge \cite{ShenZhuge} provides the optimal convergence rates also in the case $d=2,3$. Note that the latter paper also establishes convergence rates in the Neumann case which are optimal for $d\geq 3$.
%For the case of nonlinear elliptic equations or domains with less smoothness, see e.\,g.\ \cite{Arisawa,FeldmanKimSouganidis,KenigLinShen} and the references therein.

Let us recall some basic concepts and notions from homogenization of linear elliptic equations of the form \eqref{LinearEllipticEquation} on the whole space. Homogenization occurs for elliptic PDEs of the form \eqref{LinearEllipticEquation} in case of periodic or random coefficient fields $a$. In these cases, the large-scale behavior of the equation is captured by an \emph{effective equation} of the form $-\nabla \cdot (a_{hom}\nabla u_{hom})=0$, where $a_{hom}$ is a \emph{constant} effective coefficient. It is a simple observation that while affine functions $x\mapsto \xi \cdot x + c$ solve constant-coefficient equations of the form $-\nabla \cdot (a_{hom}\nabla u_{hom})=0$, they are in general not solutions to the original equation $-\nabla \cdot (a\nabla u)=0$. One may therefore attempt to ``correct'' the affine function by adding a perturbation $\phi_\xi$ which accounts for the oscillations in the coefficient field $a$ and which ensures that the resulting function $x\mapsto \xi\cdot x+c+\phi_\xi$ solves the equation $-\nabla \cdot(a\nabla u)=0$. This ansatz leads to the notion of \emph{homogenization correctors}, which by definition are solutions to the equation
\begin{subequations}
\label{CorrectorEquation}
\begin{align}
\label{CorrectorPhiEquation}
-\nabla \cdot (a\nabla \phi_\xi)=\nabla \cdot (a\xi).
\end{align}
Obviously, the corrector $\phi_\xi$ may be chosen to depend linearly on $\xi$; we shall denote the corrector corresponding to a coordinate function $x\mapsto x_i$ (i.\,e.\ to $\xi=e_i$) also by $\phi_i$. In view of the heuristics that $\xi\cdot x+\phi_\xi$ should be a perturbation of the polynomial $\xi\cdot x$, correctors are required to grow sublinearly, i.\,e.\ to satisfy a bound of the form $|\phi_\xi(x)|\ll |x|$ for large $|x|$.

The effective coefficient $a_{hom}$ is determined by the following heuristics: Consider an affine function $x\mapsto \xi\cdot x$ in the homogenized picture and the corresponding ``corrected affine function'' $x\mapsto \xi \cdot x+\phi_\xi(x)$ in the microscopic (non-homogenized) picture. Then, the flux in the homogenized picture $a_{hom}\nabla (\xi\cdot x)=a_{hom}\xi$ should correspond to the average of the flux in the microscopic (non-homogenized) picture $a\nabla (\xi\cdot x+\phi_\xi(x))=a(\xi+\nabla \phi_\xi(x))$. In stochastic homogenization, by ergodicity spatial averaging corresponds to taking the expectation. Therefore, the effective coefficient is determined by the formulas
\begin{align}
a_{hom}e_i=\mathbb{E}[a(e_i+\nabla \phi_i)].
\end{align}
Let us mention that like in periodic homogenization, the homogenized coefficient $a_{hom}$ satisfies bounds similar to \eqref{Elliptic}; see e.\,g.\ \cite{GloriaNeukammOtto}.

In quantitative homogenization, it is convenient to introduce a dual quantity to the corrector $\phi_i$, namely a vector potential $\sigma_{ijk}$ for the flux correction $a(e_i+\nabla \phi_i)-a_{hom}e_i$ (i.\,e.\ a vector potential for the difference between the flux in the microscopic picture and the flux in the homogenized picture in the case of the macroscopic affine function $x\mapsto x_i$). The vector potential $\sigma_{ijk}$ is skew-symmetric in the last two indices -- i.\,e.\ it satisfies $\sigma_{ijk}=-\sigma_{ikj}$ -- and its defining equation reads
\begin{equation}\label{def_sigma}
\nabla_{k} \cdot \sigma_{ijk} = e_j \cdot ( a (e_i + \nabla \phi_i ) -  a_{hom} e_i).
\end{equation}
\end{subequations}
Approximating a solution $u$ by the \emph{two-scale expansion}
\begin{align}
\label{TwoScaleExpansion}
u_{2-scale}:=u_{hom}+\sum_{i=1}^d \phi_i \partial_i u_{hom},
\end{align}
a simple computation shows that the error $w:=u-u_{2-scale}$ satisfies the equation
\begin{align*}
-\nabla \cdot (a\nabla w) = \nabla \cdot \left(\sum_{i=1}^d(\phi_i a - \sigma_i)\partial_i \nabla u_{hom}\right).
\end{align*}
Homogenization effects are then encoded in terms of growth estimates for the corrector $(\phi, \sigma)$: The previous formula allows one to turn estimates on the corrector $(\phi,\sigma)$ into estimates for the homogenization error. Introducing the notation
\begin{equation}\label{DefinitionDelta}
\delta_{R} := \displaystyle\frac{1}{R} \left( \displaystyle\fint_{\{|x|<R\}} |(\phi, \sigma)|^2 \,dx \right)^{1/2}
\end{equation}
as a measure for the sublinearity of the corrector at scale $R$, the qualitative convergence
\begin{align*}
\lim_{R \rightarrow \infty } \delta_R = 0
\end{align*}
is sufficient for homogenization to occur. In fact, an estimate of the form
\begin{align*}
\sup_{R\geq r^\ast}\delta_R\leq \frac{1}{C(d,\lambda,\alpha)}
\end{align*}
is sufficient for a $C^{1,\alpha}$-type excess-decay estimate on scales larger than $r^\ast$ as shown in \cite{GloriaNeukammOtto}. In the same work, the almost sure existence of correctors subject to this condition of qualitatively sublinear growth has been established, assuming only stationarity and qualitative ergodicity of the ensemble.
The higher-order regularity theory in \cite{FischerOtto} relies on the slight quantification of sublinear growth of the corrector
\begin{equation}\label{conditionFO}
\displaystyle\sum_{m = 0}^{\infty} \delta_{2^m} < \infty.
\end{equation}
To ensure the (almost sure) existence of correctors with this quantified sublinear growth, replacing the assumption of just qualitative ergodicity by a very mild assumption on decay of correlations is sufficient, see e.\,g.\ \cite{CorrectorEstimatesSlowDecorrelation}. In the case of ideal decorrelation -- e.\,g.\ finite range of dependence -- and $d\geq 3$, $\delta_r$ behaves like $\frac{1}{r}$ and therefore $\delta_{2^m}$ behaves like $2^{-m}$ , see \cite{GloriaOttoNew}.

Turning our attention to homogenization in the half-space setting, it becomes apparent that the homogeneous Dirichlet boundary conditions on $\bhs$ introduce further restrictions on the affine polynomials which are necessary to describe the behavior of solutions to the equation $-\nabla \cdot (a\nabla u)=0$ on $\hs$: Basically, the boundary conditions exclude all polynomials $\xi \cdot x$ with $\xi\cdot e_i\neq 0$ for some $i\neq d$ from playing a relevant role in the approximate description of solutions.

As a second difference to the whole-space case, correcting the remaining relevant affine polynomial $x\mapsto x_d$ with the whole-space corrector $\phi_d$ leads to a violation of the boundary conditions on $\bhs$. It becomes therefore necessary to construct a corrector $\phi_d^\mathbb{H}$ which is adapted to the half-space setting (i.\,e.\ satisfies the homogeneous Dirichlet boundary conditions on $\bhs$). In the present work, we shall present an entirely deterministic argument which modifies a given whole-space corrector $\phi_d$ to yield a corrector $\phi_d^\mathbb{H}$ that satisfies the Dirichlet boundary conditions on $\bhs$. The only condition that we need to impose on the whole-space corrector $(\phi,\sigma)$ is the (slightly stronger) condition on quantitative sublinearity
\begin{equation}\label{condition}
\sum_{m=0}^{\infty} m \delta_{2^m}^{1/3} < \infty.
\end{equation}
Note that this condition is implied by an estimate of the form
\begin{align*}
\delta_r \lesssim \frac{1}{|\log r|^{6+\epsilon}}
\quad\quad\text{for large }r
\end{align*}
for arbitrarily small $\epsilon>0$.
Again, a very mild assumption on decorrelation is sufficient to ensure the almost sure existence of whole-space correctors with this growth property \cite{CorrectorEstimatesSlowDecorrelation}. Note that as the vector potential for the flux correction $\sigma_{djk}$ is defined in terms of the corrector $\phi_d$, after modifying $\phi_d$ to obtain $\phi_d^\h$ we also need to construct an appropriately adjusted vector potential $\sigma_{djk}^\mathbb{H}$. 

{\bf Notation.} Throughout the paper, we denote the number of spatial dimensions by $d$. The notation $\hs$ is used for the half-space $\{x\in \mathbb{R}^d:x_d>0\}$. By $B_r$, we denote the ball of radius $r$ centered at the origin. The half-ball of radius $r$ centered at the origin -- i.\,e.\ the set $\{x\in \mathbb{R}^d:|x|<r,x_d>0\}$ is denoted by $B_r^+$; correspondingly, the notation $B_r^-$ is used for the set $-B_r^+$. By $B_r(x)$ we denote the ball of radius $r$ with center $x$. For two sets $M$ and $N$, the set $\{m\in M:m\notin N\}$ is denoted by $M\setminus N$.

When it is not important to keep track of constants, we use the notation ``$\lesssim$'' to mean ``up to a constant depending on $d$, $\lambda$''. The notation $C(d, \lambda, \alpha)$ is also used to denote a generic constant depending on the quantities in the brackets. By ``$a\ll b$'' we mean $a\leq \frac{1}{C(d,\lambda)} b$ for some large enough constant $C(d,\lambda)$.

For a measurable set $A\subset \mathbb{R}^d$, we denote its $d$-dimensional Lebesgue measure by $|A|$. By $\int_A f \,dx$ we denote the Lebesgue integral of the function $f$ over the set $A$. By $\fint_A f\,dx$ we denote the average integral, i.\,e.\ $\frac{1}{|A|}\int_A f\,dx$.

For a vector or tensor, the subscripts before a comma refer to components and the subscripts after a comma refer to a scale (not to taking a partial derivative): For example, $\sigma_{djk,M}^\h$ refers to the component $djk$ of a modified vector potential for the flux correction which has been adapted on scales $\leq 2^M r_0$ (with $r_0$ denoting the base scale, see Section \ref{ExpositionStrategy} below).

The function space $C^{k,\alpha}$ (with $k\in \mathbb{N}_0$ and $\alpha\in (0,1]$) consists of the functions whose derivatives up to order $k$ are (locally) H\"older continuous with exponent $\alpha$. The (possibly weak) partial derivative with respect to the $j$th coordinate will be denoted by $\partial_j$. By $\dot H^1_0(\hs)$, we denote the space of locally integrable functions $v$ with square-integrable gradient and vanishing trace on $\bhs$, equipped with the norm $||v||_{\dot H^1_0(\hs)}:=(\int_{\hs} |\nabla v|^2 \,dx)^{1/2}$.

\section{Main Results}

Our first main theorem ensures the existence of half-space-adapted homogenization correctors with the appropriate (sublinear) growth behavior. The key assumption of the theorem is the existence of correctors on the whole space which are sublinear in the mildly quantified sense (\ref{condition}).
\begin{customthm}{1}\label{existenceofcorrectors}
Let $a:\mathbb{R}^d\rightarrow \mathbb{R}^{d\times d}$ be a uniformly elliptic and bounded coefficient field in the sense \eqref{Elliptic}.
Assume that for this coefficient field $a$ there exists a whole-space corrector $(\phi, \sigma)$ satisfying the corrector equations \eqref{CorrectorEquation} and the growth condition (\ref{condition}). Then there exists a half-space-adapted corrector $(\tp, \ts)$ with the following properties:
\begin{itemize}
    \item[i)] For $i\neq d$, the correctors $\tp_i$ and $\ts_i$ coincide with the restriction of the whole-space correctors to the half-space $\phi_i|_{\hs}$ and $\sigma_i|_{\hs}$.
    \item[ii)] The corrector $\tp_d$ is adapted to the half-space setting in the sense
    \begin{subequations}
    \label{halfspacecorrector}
    \begin{align}
    \label{CorrectorEquationHalfSpace}
    ~~~~~~~~-\di (a \nabla \tp_{d}) &= \nabla \cdot (a e_d) &&\textrm{in }\quad \hs,
    \\
    ~~~~~~~~\tp_d &\equiv 0 &&\textrm{on }\quad \bhs.
    \end{align}
    \end{subequations}
    \item[iii)] $\ts_d$ is a vector potential for the flux correction corresponding to $x_d+\phi_d^\mathbb{H}$ in the sense that it is skew-symmetric and satisfies
    \begin{align}
    \label{halfspacepotential}
    \nabla_k \cdot \ts_{djk} = e_j \cdot (a(e_d + \nabla \tp_d )- a_{hom}e_d).
    \end{align}
    \item[iv)] The corrector grows sublinearly in the sense that
    \begin{align*}
    \delta^{\h}_{r} := \frac{1}{r} \left( \fint_{B_r^+} |(\phi^{\h}, \sigma^{\h}) |^2 \,dx + \fint_{B_r^-} \sum_{i=1}^{d-1} |\phi_i|^2 \,dx \right)^{1/2}
    \end{align*}
    satisfies
    \begin{align*}
    \lim_{r\rightarrow \infty} \delta^{\h}_r =0.
    \end{align*}
    In particular, for any $0<\alpha<1$ there exists a finite radius $r^\ast$ for which the condition (\ref{smallness}) below is satisfied.
\end{itemize}
\end{customthm}
In fact, our proof shows that a quantitative estimate on the sublinear growth of the whole-space corrector in the form
\begin{align*}
\delta_r \leq \frac{C}{r^\gamma}
\end{align*}
for some $\gamma \in (0,1]$ may be turned into an estimate on the half-space-adapted corrector of the form
\begin{align}
\label{QuantitativeBound}
\delta_r^\h \leq \frac{\tilde C}{r^{\gamma/3}}.
\end{align}
This bound is a consequence of more precise estimates on the right-hand sides of the inequalities \eqref{final1}, \eqref{nabla}, and \eqref{QuantitativeSublinear} in the proof below. However, one should not expect the estimate \eqref{QuantitativeBound} to be optimal, which is why we did not emphasize this quantitative bound in our theorem.

Our second main theorem transfers regularity properties from the constant-coefficient equation $-\nabla \cdot (a_{hom}\nabla u_{hom})=0$ to the equation with possibly oscillating coefficients $-\nabla \cdot (a\nabla u)=0$. The key requirement of the theorem is that approximate homogenization has occurred, as assessed by the sublinearity condition for the half-space-adapted corrector \eqref{smallness}. In this case, a large-scale regularity theory in the form of a corresponding decay estimate for the tilt-excess becomes available. As a second consequence, we infer a mean-value property for $a$-harmonic functions.
\begin{customthm}{2}\label{excessbound}
Let $a:\mathbb{R}^d\rightarrow \mathbb{R}^{d\times d}$ be a coefficient field satisfying the uniform ellipticity and boundedness assumptions \eqref{Elliptic}. For any fixed H\"{o}lder exponent $0<\alpha <1$ there exists a constant $C_{\alpha}(d, \lambda)$ such
that the following statements hold:

Suppose that for some radius $R>0$ there exist half-space-adapted homogenization correctors $(\tp, \ts)$ which satisfy the defining equations of the corrector \eqref{CorrectorEquation} on $B_R^+$, for which $\tp_d$ satisfies homogeneous Dirichlet boundary conditions on $\bhs\cap B_R$, and for which $\phi^\h_i$ for $i\neq d$ is the restriction of a corrector $\phi_i$ on $B_R$ to $B_R^+$. Suppose furthermore that the correctors $(\tp,\ts)$ are sublinear on larger scales in the sense that the quantity
\begin{align*}
\delta^{\h}_{r} := \frac{1}{r} \left(\fint_{B_r^+} |(\phi^{\h}, \sigma^{\h}) |^2 \,dx +\fint_{B_r^-} \sum_{i=1}^{d-1} |\phi_i|^2 \,dx \right)^{1/2}
\end{align*}
satisfies an estimate of the form
\begin{align}
\label{smallness}
\delta^{\h}_{r} \leq \frac{1}{C_{\alpha} (d, \lambda)} \textit{ for all $r\geq r^*$}
\end{align}
for some radius $0<r^*<R$.

Let $u\in H^1(B_R^+)$ be an $a$-harmonic function with homogeneous Dirichlet boundary conditions on $\bhs$, i.\,e.\ let $u$ be a solution to the problem
\begin{align*}
-\di (a \nabla u) &= 0 \quad\quad\quad \textrm{ in }\quad B_R^+,
\\
u &\equiv 0 \quad\quad\quad \textrm{ on }\quad \partial \hs \cap B_R.
\end{align*}
Introduce the half-space-adapted tilt-excess
\begin{align*}
\textrm{Exc}^{\h}(r):=\inf_{b\in \mathbb{R}} \fint_{B_r^+} |\nabla u-b(e_d+\nabla\tp_d)|^2 \,dx.
\end{align*}
Then for any $r\in [r^\ast,R]$ the excess-decay estimate
\begin{align}
\label{ExcessDecay}
 \textrm{Exc}^{\h}(r) \lesssim \left( \displaystyle\frac{r}{R} \right)^{2\alpha} \textrm{Exc}^{\h}(R)
\end{align}
is satisfied.

Furthermore, for $r\in [r^\ast,R]$ the mean-value property
\begin{align}
\label{mvp}
 \firs{r}{\nabla u} \leq C_{Mean}(d, \lambda) \firs{R}{\nabla u}
\end{align}
holds for some constant $C_{Mean}(d,\lambda)$ depending only on the dimension $d$ and the ellipticity constant $\lambda$.

Finally, for all $r\in [r^\ast,R]$ the tilt-excess functional
\begin{align*}
\fint_{B_r^+} |\nabla u - b(e_d +\nabla \tp_{d})|^2 \,dx
\end{align*}
is coercive as a function of $b \in \mathbb{R}$ in the sense
\begin{align}
\label{Coercivity}
\fint_{B_r^+} |\nabla u - b(e_d +\nabla \tp_{d})|^2 \,dx \geq c(d,\lambda) |b-b_{min}|^2
\end{align}
for some $b_{min}\in \mathbb{R}$.
\end{customthm}
Combining Theorem \ref{existenceofcorrectors} with Theorem \ref{excessbound} yields the following Liouville principle:
\begin{customcorollary}{1.1}\label{Liouville}
Let $a:\mathbb{R}^d\rightarrow \mathbb{R}^{d\times d}$ be a coefficient field which is uniformly elliptic and bounded in the sense \eqref{Elliptic}. Suppose that for the coefficient field $a$ homogenization correctors $(\phi,\sigma)$ exists which satisfy the corrector equations \eqref{CorrectorEquation} and the growth condition $(\ref{condition})$. Then, there exists a sublinearly growing homogenization corrector on the half-space $\phi_d^\h$ in the sense that it satisfies \eqref{halfspacecorrector} and
\begin{align*}
\lim_{r\rightarrow \infty} \frac{1}{r}\bigg(\fint_{\{|x|<r\}}|\phi_d^\h|^2 \,dx\bigg)^{1/2}
=0.
\end{align*}
Furthermore, any $a$-harmonic function $u\in H^1_{loc}(\hs)$ with homogeneous Dirichlet boundary conditions on $\bhs$ and subquadratic growth in the sense
\begin{align}
\label{subquad}
\lim_{r \rightarrow \infty}\frac{1}{r^{1+\alpha}}\bigg(\firs{r}{u}\bigg)^{1/2}=0
\end{align}
for some $0<\alpha<1$ must be of the form
\begin{align*}
u = b \cdot (x_d + \tp_d)
\end{align*}
for some $b\in \mathbb{R}$.
\end{customcorollary}

\subsection{Strategy for the construction of half-space-adapted correctors}
\label{ExpositionStrategy}
In the present section we give an exposition of our strategy for the construction of half-space-adapted homogenization correctors (Theorem \ref{existenceofcorrectors}). At several points, it will be important to keep track of certain constants in the estimates:
\begin{itemize}
 \item $C_{Mean}(d,\lambda)$, which comes from the mean-value property (\ref{mvp}),
 \item $C_P(d)$, which we take to be an upper bound for the Poincar\'{e} constant of the unit ball in $\mathbb{R}^d$ with homogeneous Dirichlet boundary conditions and the Poincar\'{e} constant of the unit half-ball $B_1^+$ with homogeneous Dirichlet boundary conditions on $\bhs\cap B_1$,
 \item and $C_I(d)$, which comes from the constant-coefficient regularity estimate (\ref{s(5)}) below.
\end{itemize}
We assume that all of these constants are larger than $1$.

\smallskip\smallskip
\noindent
\textbf{Step 1: Construction of a sublinear $\tp_d$ up to a certain scale.}
\smallskip
\\Our approach for the construction of the half-space-adapted corrector $\tp_d$ is to adapt the whole-space corrector $\phi_d$ to the Dirichlet boundary conditions on $\bhs$. We would like to achieve this by subtracting from $\phi_d$ a sublinearly growing function $\tilde{\varphi}$ that is $a$-harmonic on $\hs$ and equals $\phi_d$ on the boundary, i.\ e.\ by setting $\tp_d := \phi_d - \tilde{\varphi}$ with $\tilde{\varphi}$ being a sublinearly growing solution to the problem
\begin{subequations}
\label{TildeVarphiProblem}
\begin{align}
\label{p(1)a}
~~~~~~~~~~~-\di (a \nabla \tilde{\varphi}) &= 0 && \textrm{in } \hs,
\\
\label{p(1)b}
\tilde{\varphi} &= \phi_d &&  \textrm{on } \bhs.
\end{align}
\end{subequations}

As (\ref{p(1)a}) is a linear equation, we can decompose the right-hand side in \eqref{p(1)b} into contributions from dyadic annuli, solve the corresponding problems, and then add the solutions to obtain $\tilde \varphi$. We will show that this sum converges and sums to a sublinearly growing function.

Pursuing this strategy, let $r_0=2^{m_0}$, $m_0\in \mathbb{N}$, be a generic dyadic radius. Let $\{ \eta_m | -1 \leq m \}$ be a radial partition of unity with $\operatorname{supp} \eta_{-1}\subset \{x\in \mathbb{R}^{d}:|x|\leq r_0\}$ and $\operatorname{supp} \eta_m \subset \{x\in \mathbb{R}^{d}:r_0 2^{m-1}\leq |x|\leq r_0 2^{m+1}\}$ for $m\geq 0$; suppose that $\eta_m$ satisfies an estimate of the form $|\nabla \eta_m|\leq \frac{4}{r_0 2^m}$. Also, for $L_m \in (0,r_0 2^{m+1}]$ consider one-dimensional cutoff functions $S_m(x)=S_m(x_d)$ satisfying $S_m(x)=1$ for $|x_d|\leq L_m$ and $S_m(x)=0$ for $|x_d|\geq 2L_m$; suppose that $|\nabla S_m|\leq \frac{2}{L_m}$. Note that we shall later choose $L_m\ll r_0 2^{m+1}$.

Introducing the cutoffs $\chi_m (x) := \eta_m(x) S_m(x)$, we then consider the Lax-Milgram solutions $\varphi_m\in H^1_0(\hs)$ to the problem
\begin{subequations}
\label{p(2)}
\begin{align}
\label{p(2)a}
~~~~~~~~ - \di (a \nabla \varphi_m) &= \di (a \nabla(\chi_m \phi_d )) && \textrm{ in } \quad \hs,
\\
\varphi_m &= 0 && \textrm{ on }\quad \bhs.
\end{align}
\end{subequations}
Defining $\varphi_M^{\Sigma} := \sum_{m=-1}^M \varphi_m $ and $\tilde{\varphi}_M^{\Sigma}  := \varphi_M^{\Sigma} + \sum_{m=-1}^M \chi_m \phi_d$, we see that
\begin{align*}
\phi^{\h}_{d,M} := \phi_d - \tilde{\varphi}^{\Sigma}_M
\end{align*}
solves the corrector equation \eqref{CorrectorEquationHalfSpace} in $\hs$ with homogeneous Dirichlet boundary conditions on $\bhs \cap B_{r_0 2^M}$.

In order to estimate the size of the modification $\tilde{\varphi}_M^\Sigma$ on a half-ball $B_r^+$, we will first deduce an estimate for the ``near-field contributions'', i.\,e.\ the $\varphi_m$ for which the inclusion $\operatorname{supp} \chi_m \subseteq B_{16r}$ holds. As we shall see, this is easily done with the standard energy estimate for the equation (\ref{p(2)}) and an appropriate estimate for $\chi_m \phi_d$. The energy norm of the term $\chi_m \phi_d$ in turn may be made small by an appropriate choice of $L_m$.

\begin{customlemma}{2.1}
\label{first}
Let the assumptions of Theorem \ref{existenceofcorrectors} be satisfied. Let $m\geq -1$. Then there exists $L_m\ll r_0 2^{m+1}$ and a constant $C_1(d,\lambda)$ such that the following is true: For any $r>0$ the estimates
\begin{align}
\label{chiphiformula}
\fir{r}{\nabla( \chi_m \phi_d)} \leq C_1(d, \lambda) \left(\frac{r_0 2^{m+1}}{r} \right)^{d/2} \delta_{r_0 2^{m+1}}^{1/3}
\end{align}
and
\begin{align}
\label{energy1formula}
 \fir{r}{\nabla \varphi_m} \leq C_1(d, \lambda) \left( \frac{r_0 2^{m+1}}{r} \right)^{d/2} \delta_{r_0 2^{m+1}}^{1/3}
\end{align}
hold. In particular, for any $r\geq \frac{1}{16} r_0 2^{m+1}$ the function $\varphi_m$ satisfies the bound
\begin{align}
\label{assumption1}
\fir{r}{\nabla \varphi_m} \leq C_2(d, \lambda)
\min \bigg\{1, \left(\frac{r_0 2^{m+1}}{r}\right)^{d/2} \bigg\} \delta_{r_0 2^{m+1}}^{1/3}
\end{align}
with $C_2 := C_{Mean} C_1 8^{d}$.
\end{customlemma}
%In particular, this yields that $\phi^{\h}_{d,M}$ as defined above is sublinear when $ \frac{r_0 2^{M+1}}{r} \leq 16$.
However, we will need the estimate \eqref{assumption1} on $B_r^+$ also for the ``far-field contributions'', i.\,e.\ for the $\varphi_m$ for which $\operatorname{supp} \chi_m \cap B_{4r}^+ = \emptyset$ holds. For such $m$ with $r_0 2^{m+1}\geq 16r$, the estimate \eqref{assumption1} will be established in Step 3 below.

\smallskip\smallskip
\noindent
\textbf{Step 2: Construction of a sublinearly growing $\ts_d$ up to a certain scale.}
\smallskip
\\Having constructed a corrector $\phi_{d,M}^\h$ which satisfies the homogeneous Dirichlet boundary conditions on $\bhs\cap B_{2^{M}r_0}$, we need to construct a corresponding vector potential $\ts_{d,M}$ for the flux correction, as the vector potential for the flux correction depends on the corrector through its defining equation \eqref{def_sigma}. Again, our approach is to adapt the vector potential $\sigma_d$ to take into account the modification $\phi_{d,M}^\h-\phi_d$ of the corrector by adding a correction $\psi_{jk,M}$: We construct sublinearly growing functions $\psi_{jk,M}$ that satisfy
\begin{align}
\label{equationforpsi}  
    - \nabla_k \cdot \psi_{jk,M} = e_j \cdot \big(a (e_d+\nabla \phi_{d,M}^\h) - a (e_d+\nabla \phi_d) \big) \quad \textrm{ in } \hs 
\end{align}
and define
\begin{align*}
\ts_{djk,M} := \s_{djk} - \psi_{jk,M}.
\end{align*}
Note that in order to ensure the skew-symmetry of $\ts_{d,M}$, we need to construct the $(\psi_{jk,M})_{jk}$ as skew-symmetric.
It turns out that a suitable ansatz is
\begin{align}
\label{AnsatzPsi}
\psi_{jk,M} := \partial_k v_{j,M}-\partial_j v_{k,M}
\end{align}
with $v_{,M}: \hs \rightarrow \mathbb{R}^d$ solving the equation
\begin{subequations}
\label{s(3)}
\begin{align}
\label{s(3)a}
-\Delta v_{j,M} &=  e_j\cdot (a (e_d+\nabla \phi_{d,M}^\h) - a (e_d+\nabla \phi_d)) \quad &&\textrm{in } \hs,
\\
v_{j,M} &= 0 \quad &&\textrm{for $j \neq d$  on } \bhs,
\\
\partial_d v_{d,M} &= 0 \quad &&\textrm{on } \bhs.
\end{align}
\end{subequations}
First, note that the skew-symmetry of $\psi_{jk,M}$ is built into the ansatz \eqref{AnsatzPsi}. Furthermore, differentiating the equation (\ref{s(3)}), we infer
\begin{subequations}
\label{s(6)}
\begin{align}
- \Delta (\nabla_k \cdot v_{k,M} ) &= 0 && \textrm{ in } \hs,
\\
\nabla_k \cdot v_{k,M}  &= 0 && \textrm{ on } \bhs.
\end{align}
\end{subequations}
By the Liouville principle for harmonic functions with homogeneous Dirichlet boundary conditions on $\hs$, sublinear growth of $\nabla_k \cdot v_{k,M} $ entails that $\nabla_k \cdot v_{k,M}   \equiv 0$. This leads, as desired, to the conclusion
\begin{align*}
-\nabla_k \cdot \psi_{jk,M} &=  \sum_{k=1}^d (\pk \pj v_{k,M} - \pk^2 v_{j,M})
\\&
= \pj (\nabla_k \cdot v_{k,M})  - \Delta v_{j,M}
\\&
= e_j\cdot (a (e_d+\nabla \phi_{d,M}^\h) - a (e_d+\nabla \phi_d)).
\end{align*}
To summarize, in order to obtain a solution to \eqref{equationforpsi} it suffices to construct solutions $v_{j,M}$ to (\ref{s(3)}) for which $\nabla_k \cdot v_{k,M}$ is a sublinearly growing function (note that we shall actually prove the stronger statement of sublinear growth of $\nabla v_{k,M}$).

To construct such a solution $v_{j,M}$, notice that, as $\phi_{d,M}^{\h} - \phi_d$ is $a$-harmonic on $\hs$, we may rewrite the right-hand side in \eqref{s(3)a} as
\begin{align*}
e_j\cdot (a (\nabla \phi_{d,M}^{\h} - \nabla \phi_d)) &=
e_j \cdot a (\nabla \phi_{d,M}^{\h} - \nabla \phi_d) + x_j \di (a (\nabla \phi_{d,M}^{\h} - \nabla \phi_d))
\\&
=\di (x_j a (\nabla \phi_{d,M}^{\h} - \nabla \phi_d)).
\end{align*}
Our strategy, just like in Step 1, is now to work with a decomposition into contributions from dyadic annuli: Reusing the partition of unity $\eta_m$ from Step 1, we consider the Lax-Milgram solutions $v_{j,M}^n$ of the problems
\begin{subequations}
\label{s(2)}
\begin{align}
- \Delta v_{j,M}^n &= \di ( \eta_n x_j a ( \nabla \phi^{\h}_{d,M} - \nabla \phi_d )) && \textrm{in } \hs,
\\
v_{j,M}^n &= 0 && \textrm{for } j \neq d \textrm{ on } \bhs,
\\
\partial_d v_{d,M}^n(x) &= 0 &&  \textrm{on } \bhs.
\end{align}
\end{subequations}
Here, in order to find the solutions $v_{j,M}^n$ for $j\neq d$ we apply Lax-Milgram to the space ${\dot H^1_0}(\hs)$. To find the solution $v_{d,M}^n$, we apply Lax-Milgram to the space of locally integrable functions $v$ with square-integrable gradient subject to the constraint $\fint_{B_{r_0}^+} v \,dx=0$; we equip this space with the norm $||v||:=(\int_{\hs} |\nabla v|^2\,dx)^{1/2}$.

In order to obtain $v_{j,M}$, we intend to sum all of the contributions. However, to ensure that on a half-ball $B_r^+$ the ``far-field contributions'' -- i.\,e.\ the $v_{j,M}^n$ with $2^{n+1} r_0\geq 16r$ -- do not destroy the smallness of the sum $\sum_{n=-1}^\infty \nabla v_{j,M}^n$, we must enforce ``quadratic'' behavior of $v_{j,M}^n$ around the origin by subtracting off the linear growth of $v_{j,M}^n$: Set
%enforce the sublinear of $v_{j,M}$ we subtract-off the initial linear growth of the $v_{,M}^n$ and, therefore, define
\begin{align}
\label{DefinitionbjMn}
b_{j,M}^n = :\bigg\{
\begin{array}{ll}
       0 & \textrm{  if  } n=-1 \\
     \nabla v_{j,M}^n(0) & \textrm{  if  } n \neq -1
\end{array}.
\end{align}
Notice that $b_{jk,M}^n = 0$ unless $n \neq -1$ and either $j=d$ and $k\neq d$ or $j \neq d$ and $k=d$.
We then obtain the following estimate which in particular shows that $v_{j,M}^n-b_{j,M}^n\cdot x$ indeed displays quadratic behavior in the interior $\{|x|<2^n r_0\}$:
\begin{customlemma}{2.2}
\label{2.2}
Let the assumptions of Theorem \ref{existenceofcorrectors} be satisfied.
Let $M \geq -1$ and $n\geq -1$. Then for any $r \geq r_0$ and any $j,k \in \left \{1,..., d \right\}$ we have the estimate
\begin{align*}
&\df \fir{r}{\partial_k ( v_{j,M}^n - b_{j,M}^n \cdot x)}\\
&\leq C_3(d, \lambda )\min\bigg\{ 1, \frac{r_0 2^{n+1}}{r}\bigg\} \fir{r_0 2^{n+1}}{\cd}
\end{align*}
with $C_3(d, \lambda) := 4 C_4 C_I$.
\end{customlemma}
This estimate immediately enables us to pass to the limit $N\rightarrow \infty$ in the sum $\sum_{n=-1}^N (v_{j,M}^n - b_{j,M}^n \cdot x)$.
\begin{customlemma}{2.3}
\label{2.3}
Let the assumptions of Theorem \ref{existenceofcorrectors} be satisfied and let the $L_m$ be chosen as in Lemma \ref{first}. Then for $r\geq r_0$ and $j \in \left \{1,..., d \right\}$ the series $\sum_{n=-1}^\infty (v_{j,M}^n - b_{j,M}^n \cdot x)$ converges absolutely in $H^1(B_r^+)$ to a limit $v_{j,M}$. For this limit, the function $\psi_{jk,M}=\partial_k v_{j,M}-\partial_j v_{k,M}$ satisfies the equation
\begin{align}
\label{EquationPsijkM}
 - \nabla_k \cdot \psi_{jk,M} = e_j\cdot a(\nabla \phi_{d,M}^\h-\nabla \phi_d)  \quad\quad\textrm{ in } \hs
\end{align}
and for any $r \geq r_0$ and any $j,k \in \left\{ 1,..., d \right\}$ we have the estimate
\begin{align}
\nonumber
\df \fir{r}{\psi_{jk,M}}
\leq 2 C_3 (d, \lambda) \sum_{n=-1}^{\infty}&\min\bigg\{1, \frac{r_0 2^{n+1}}{r} \bigg\}
\\&~~\times
\fir{r_0 2^{n+1}}{\cd}.
\label{result1}
\end{align}
\end{customlemma}

\smallskip\smallskip
\noindent\textbf{Step 3: Inductively building a sublinear corrector on larger scales.}
\smallskip
\\
Notice that in the previous two steps the radius $r_0$ was arbitrary. In the present step, we now choose $r_0$ independently of $m$ in such a way that the estimate (\ref{assumption1}) does not only hold for $r\geq \frac{1}{16}r_0 2^{m+1}$, but more generally for arbitrary $r\geq r_0$.

To extend the inequality \eqref{assumption1} for $\varphi_{m+1}$ to arbitrary $r\geq r_0$, we shall crucially rely on the mean-value property \eqref{mvp} for $a$-harmonic functions for radii $r\in [r_0,r_0 2^m]$. To this aim, we proceed by induction in $m$; to show \eqref{assumption1} for $\varphi_{m+1}$, we shall use the already-constructed corrector $(\tp_{d,m}, \ts_{d,m})$ and establish that it satisfies the estimate \eqref{smallness}, which by Theorem \ref{excessbound} entails the mean-value property \eqref{mvp} for $a$-harmonic functions on scales $r\in [r_0,r_0 2^m]$ with $R=r_0 2^m$.

We therefore have to choose the dyadic radius $r_0=2^{m_0}$ in such a way that we obtain a bound which guarantees for all $m$ that the smallness condition \eqref{smallness} is satisfied by $(\tp_{d,m}, \ts_{d,m})$.

\begin{customlemma}{2.4}
\label{induct}
Let the assumptions of Theorem \ref{existenceofcorrectors} be satisfied -- in particular, suppose that for the coefficient field $a$ there exist whole-space correctors which satisfy the quantitative sublinear growth condition \eqref{condition} -- and let the $L_m$ be chosen as in Lemma~\ref{first}.
Then there exists $r_0>0$ independent of $M\in \{-1,0,1,2,\ldots\}$ with the following property: If the $\varphi_m$ satisfy the estimate
\begin{align}
\label{assumption_new}
\fir{r}{ \nabla \varphi_m} \leq C_2
\min\bigg\{1, \left(\frac{r_0 2^{m+1}}{r}\right)^{d/2}\bigg\}\delta_{r_0 2^{m+1}}^{1/3}
\end{align}
for all $r \geq r_0$ and all $m\in \{-1,\ldots,M\}$ (recall the definition $C_2 := C_{Mean} C_1 8^{d}$), then $(\tp_{d,M}, \ts_{d,M})$ satisfies the smallness condition (\ref{smallness}) for $\alpha = 1/2$ and all $r \geq r_0$, i.\,e.\ we have
\begin{align*}
\delta_{r}^{\h} \leq \displaystyle\frac{1}{C_\frac{1}{2}(d, \lambda)}.
\end{align*}
As a consequence, in this case $\varphi_{M+1}$ also satisfies the estimate (\ref{assumption_new}) for all $r\geq r_0$.
\end{customlemma}
Note that the start of the induction -- i.\,e.\ the estimate \eqref{assumption_new} for $m=-1$ -- is provided by Lemma \ref{first}.

\smallskip\smallskip
\noindent \textbf{Step 4: Passage to the limit in $M$.}
\smallskip
\\
In the last step, we pass to the limit $M\rightarrow \infty$ to obtain the half-space-adapted correctors $\phi_{d}^{\h}$ and $\sigma_{d}^{\h}$ as the limits of the sequences $\phi_{d,M}^{\h}$ and $\sigma_{d,M}^{\h}$, thereby establishing Theorem \ref{existenceofcorrectors}.

\section{Adaption of the Correctors to the Half-Space Setting}

\subsection{Step 1 -- Estimates for the modification of the corrector $\phi_d$ in the near-field case}

Lemma \ref{first} is basically a consequence of appropriate energy estimates for the defining equation of $\varphi_m$ and a suitable bound for $\chi_m \phi_d$.
\begin{proof}[Proof of Lemma \ref{first}]
Let us abbreviate $R:=r_0 2^{m+1}$.
Testing (\ref{p(2)a}) with $\varphi_m$, making use of the fact that $\varphi_m$ vanishes on $\bhs$, and estimating using the uniform ellipticity and boundedness of $a$ yields
\begin{align}
\label{energy}
\irtr{\hs}{\nabla \varphi_m} \lesssim
\ir{R}{\phi_d \nabla \chi_m} + \ir{R}{\chi_m \nabla \phi_d}.
\end{align}
We treat the two terms  on the right hand side separately. For the first, using our definition of $\chi_m$ and $L_m \leq R$, we find that
\begin{align}
\label{p(6)}
\ir{R}{\phi_d \nabla \chi_m} & \lesssim& \displaystyle\frac{R^{d/2}}{L_m} \fir{R}{\phi_d}
\leq \displaystyle\frac{R^{\frac{d+2}{2}}}{L_m} \delta_R.
\end{align}
Let us now even-reflect $\chi_m$ such that it is defined on $\mathbb{R}^d$. We may then test the corrector equation \eqref{CorrectorPhiEquation} with $\chi_m^2 (\phi_d + x_d)$. After using Young's inequality and the uniform ellipticity of $a$, this yields 
\begin{align}
\label{p(3)}
\int_{\mathbb{R}^d} \chi_m^2 |\nabla \phi_d + e_d|^2 \,dx
\lesssim \int_{\mathbb{R}^d} |\nabla \chi_m|^2 |\phi_d +x_d|^2 \,dx.
\end{align}
Now notice that we have $\operatorname{supp} \chi_m \subset [-R,R]^{d-1}\times [-2L_m,2L_m]$; in particular, on $\operatorname{supp} \chi_m$ we have $|x_d|\leq 2L_m$. The triangle inequality in $L^2(B_R)$ along with Young's inequality, the estimate (\ref{p(3)}), and the bound $|\nabla \chi_m|\leq \frac{C}{L_m}$ then yield
\begin{align*}
&\int_{B_R} |\chi_m \nabla \phi_d|^2 \,dx
\lesssim
\int_{B_R} \chi_m^2 \,dx+\int_{B_R} \chi_m^2 |\nabla \phi_d+e_d|^2 \,dx
\\&
\overset{(\ref{p(3)})}{\lesssim} |\operatorname{supp} \chi_m| + \frac{1}{L_m^2} \int_{\operatorname{supp} \chi_m} |\phi_d|^2 +|x_d|^2 \,dx
\\&
\lesssim |\operatorname{supp} \chi_m| + \frac{R^d}{L_m^2} \displaystyle\fint_{B_R}|\phi_d|^2 \,dx
\\&
\lesssim R^{d-1}L_m + \frac{R^{d+2}}{L_m^2} \delta_R^2.
\end{align*}
The second term on the right-hand side of (\ref{energy}) is therefore estimated by
\begin{align}
\label{p(5)}
\ir{R}{\chi_m \nabla \phi_d} \lesssim R^{(d-1)/2} L_m^{1/2} + \frac{R^{(d+2)/2}}{L_m}\delta_R.
\end{align}
Together, (\ref{p(5)}), (\ref{p(6)}), and  (\ref{energy}) give that
\begin{align*}
&\fir{r}{\nabla \varphi_m}
 + \fir{r}{\nabla( \chi_m \phi_d)} 
\\&
 \lesssim \left( \frac{R}{r} \right)^{d/2} \frac{R}{L_m} \delta_{R} 
 + \left( \frac{R}{r}\right)^{d/2} \left( \frac{L_m}{R} \right)^{1/2}.
\end{align*}
Choosing $L_m:= \epsilon R=\epsilon r_0 2^{m+1}$, we can optimize this expression in $\epsilon$. Plugging in the optimal $\epsilon = \delta_{R}^{2/3}$ yields
\begin{align*}
\fir{r}{\nabla \varphi_m}
+ \fir{r}{\nabla( \chi_m \phi_d)} 
&\leq C_1 \left(\frac{r_0 2^{m+1}}{r}\right)^{d/2} \delta_{r_0 2^{m+1}}^{1/3}.
\end{align*}
This directly gives \eqref{chiphiformula} and \eqref{energy1formula}. By the definition of $C_2$, for $r\geq \frac{1}{16}r_0 2^{m+1}$ this also entails the estimate \eqref{assumption1}.
\end{proof}

\subsection{Step 2 -- Estimates for the modification of the vector potential $\sigma$}
The following bound forms the basis for the estimates on the size of the modification $\psi_{jk}$ of the flux correction $\sigma_d$. It is obtained by an energy estimate for $v_{j,M}^n$ and a mean-value property of harmonic functions.
\begin{customlemma}{3.2}
\label{energy2}
Using the notation from Section \ref{ExpositionStrategy}, let $M \geq -1$, $n \geq -1$, and abbreviate $R := r_0 2^{n+1}$. Then there exists a constant $C_4=C_4(d)$ such that for any $r\geq \frac{1}{16} R$ the estimate
\begin{align*}
\fir{r}{\nabla v_{j,M}^n - b_{j,M}^n} \leq C_4 R \fir{R}{\cd}
\end{align*}
holds.
\end{customlemma}
\begin{proof}
Notice that the weak formulation of equation (\ref{s(2)}) reads
\begin{align*}
\int_{\hs} \nabla v_{j,M}^n \cdot \nabla w \,dx
= -\int_{\hs} \eta_n x_j a ( \nabla \phi^{\h}_{d,M} - \nabla \phi_d ) \cdot \nabla w \,dx
\end{align*}
for any test function $w\in H^1_0(\hs)$ in case $j\neq d$ respectively any $w\in H^1(\hs)$ in case $j=d$. In this weak formulation, no boundary terms appear: For $j \neq d$, this is a consequence of the homogeneous Dirichlet boundary conditions satisfied by the test functions on $\bhs$. For $j=d$, this is a consequence of the homogeneous Neumann boundary condition $\partial_d v_{d,M}^n =0$ on $\bhs$ and the fact that $x_d=0$ on $\bhs$.
Testing this weak formulation with $v_{j,M}^n$ and using the property $\operatorname{supp} \eta_n \subset \{|x|\leq R\}$ of the cutoff $\eta_n$ as well as the boundedness of $a$ (see \eqref{Elliptic}), we obtain the energy estimate
\begin{align}
\left(\irt{\hs}{\nabla v_{j,M}^n}\right)^{1/2}
\leq R \ir{R}{\cd}.
\label{cauchyref}
\end{align}
Using the fact that for $n\neq -1$ the functions $\partial_k v_{j,M}^n$ are harmonic in $\{|x|<\frac{R}{4}\}$ with homogeneous Dirichlet or Neumann boundary conditions on $\bhs\cap \{|x|<\frac{R}{4}\}$ (depending on $j$ and $k$) and therefore satisfy a mean-value property, we deduce by \eqref{DefinitionbjMn}
\begin{align*}
|b_{j,M}^n|\leq |\nabla v_{j,M}^n(0)|&\leq C(d) \left(\fint_{B_{R/4}} |\nabla v_{j,M}^n|^2 \,dx\right)^{1/2}
\\&
\leq C(d) R \fir{R}{\cd}.
\end{align*}
The lemma is now an easy consequence of these two estimates.
\end{proof}

Our next goal is to prove Lemma \ref{2.2}. To this aim, recall the following basic fact about harmonic functions: For any harmonic function $w$ on $B_R^+$ with either homogeneous Dirichlet or homogeneous Neumann boundary conditions on $\bhs\cap B_R$, for any $r \in (0,R/4]$ we have
\begin{align}
\label{s(5)}
\fir{r}{w - w(0)} \leq   C_I(d) \frac{r}{R} \fir{R}{w}.
\end{align}
This inequality follows from the regularity estimate (\ref{innerregularityhb}) below and the Caccioppoli estimate for harmonic function on $B_R^+$ with homogeneous Neumann or Dirichlet boundary conditions on $\bhs\cap B_R$ (for the Dirichlet case, see Lemma~\ref{caclemma}; the proof in the Neumann case is completely analogous).

\begin{proof}[Proof of Lemma \ref{2.2}]
For a given radius $r$, we separately consider the case of a ``near-field contribution'' -- defined as contributions for which $n$ satisfies $r_02^{n+1} \leq 16r$ -- and the case of a ``far-field contribution'', i.\,e.\ the case $r_02^{n+1}>16r$. Notice that, since $r \geq r_0$, $n=-1$ always corresponds to a near-field contribution.

For the near-field contributions, by Lemma \ref{energy2} we have the estimate
\begin{align}
\nonumber
&\df \fir{r}{\pk(v_{j,M}^n- b_{j,M}^n \cdot x)}
\\& \nonumber
\leq C_4\frac{r_0 2^{n+1}}{r} \fir{r_0 2^{n+1}}{\cd}
\\&
\leq 16 C_4 \min\bigg\{1,\frac{r_0 2^{n+1}}{r} \bigg\}\fir{r_0 2^{n+1}}{\cd}.
\label{nearfield1}
\end{align}
Next we address the far-field contributions, i.\,e.\ the contributions with $\frac{r_02^{n+1}}{r} > 16$. Notice that $\partial_k v_{j,M}^n - b_{jk,M}^n$ is harmonic in $B_{r_02^{n-1}}^+$ and satisfies either homogeneous Dirichlet or homogeneous Neumann boundary conditions on $\bhs\cap B_{r_02^{n-1}}$ (depending on $j$ and $k$). Furthermore, we have $\partial_k v_{j,M}^n(0)- b_{jk,M}^n=0$ and $r \leq r_0 2^{n-3}$. Therefore, an application of (\ref{s(5)}) to $w := \pk v_{j,M}^n - b_{jk,M}^n$ followed by Lemma \ref{energy2} -- the latter applied with $r:=r_0 2^{n-1}$ and $R:=r_0 2^{n+1}$ -- yields the desired bound
\begin{align}
\nonumber
&\df \fir{r}{\pk(v_{j,M}^n-b_{j,M}^n \cdot x)}
\\& \nonumber
\leq C_I \frac{1}{r_0 2^{n-1}} \fir{r_0 2^{n-1}}{\pk v_{j,M}^n - b_{jk,M}^n}
\\& \nonumber
\leq 4 C_4 C_I \fir{r_0 2^{n+1}}{\cd}.
%\label{s(7)}
\end{align}
\end{proof}

\begin{proof}[Proof of Lemma \ref{2.3}]
By Lemma \ref{2.2}, for any $r>0$ absolute convergence in $H^1(B_r^+)$ of the series
\begin{align*}
\sum_{n=-1}^\infty (v_{j,M}^n-b_{j,M}^n\cdot x)
\end{align*}
towards a limit $v_{j,M}$ follows once we have established an estimate of the form
\begin{align}
\label{SumSigmaFinite}
\sum_{n=-1}^\infty \fir{r_0 2^{n+1}}{\cd}<\infty.
\end{align}
Note that since $v_{j,M}^n$ is a weak solution of \eqref{s(2)}, the difference $v_{j,M}^n-b_{j,M}^n\cdot x$ is also a weak solution of \eqref{s(2)}. One may therefore pass to the infinite sum in the weak formulation of the problems \eqref{s(2)} (with $v_{j,M}^n$ replaced by $v_{j,M}^n-b_{j,M}^n\cdot x$) to conclude that the limit $v_{j,M}$ is a weak solution of the equation \eqref{s(3)}. Here, as test functions one uses smooth functions with bounded support in $\mathbb{R}^d$ (case $j=d$) respectively with compact support in $\hs$ (case $j\neq d$).

Lemma \ref{2.2} also implies the bound
\begin{align*}
&\df \fir{r}{\partial_k v_{j,M}}
\\
&\leq C_3(d, \lambda )\sum_{n=-1}^\infty \min\bigg\{ 1, \frac{r_0 2^{n+1}}{r}\bigg\} \fir{r_0 2^{n+1}}{\cd}.
\end{align*}
Thus, the estimate \eqref{result1} is a direct consequence of Lemma \ref{2.2}.
Furthermore, once we have established an estimate of the form \eqref{SumSigmaFinite}, this bound also entails sublinear growth of the function $\nabla_k \cdot v_{k,M}$ in the sense
\begin{align*}
\lim_{r\rightarrow \infty} \df \fir{r}{\nabla_k \cdot v_{k,M}} =0.
\end{align*}
Recalling the derivation of \eqref{equationforpsi} in the discussion of Step 2 in Section \ref{ExpositionStrategy}, we then deduce that $\psi_{jk,M}$ indeed satisfies \eqref{EquationPsijkM}.

It therefore only remains to show \eqref{SumSigmaFinite}. For any $m\in \{-1,\ldots,M\}$, the bounds \eqref{chiphiformula} and \eqref{energy1formula} -- applied with $r:=r_0 2^{n+1}$ -- entail that
\begin{align*}
\left(\fint_{B_{r_0 2^{n+1}}^+} |\nabla \varphi_m|^2 \,dx\right)^{1/2}
+\left(\fint_{B_{r_0 2^{n+1}}^+} |\nabla (\chi_m \phi_d)|^2 \,dx\right)^{1/2}
\lesssim 2^{d(m-n)/2} \delta_{r_0 2^{m+1}}^{1/3}.
\end{align*}
Taking the sum with respect to $m$ and recalling that $\phi_{d,M}^\h-\phi_d=-\sum_{m=-1}^M (\varphi_m+\chi_m \phi_d)$, we get
\begin{align*}
\fir{r_0 2^{n+1}}{\cd}
\lesssim 2^{-dn/2} \sum_{m=-1}^M 2^{dm/2} \delta_{r_0 2^{m+1}}^{1/3}.
\end{align*}
This directly implies \eqref{SumSigmaFinite}.
\end{proof}

\subsection{Step 3 -- Estimates for the modification of the corrector $\phi_d$ in the far-field case}

\begin{proof}[Proof of Lemma \ref{induct}]
For the moment, let $r_0=2^{m_0}>0$ be an arbitrary dyadic radius for which the $\varphi_{m}$ with $m\in \{-1,\ldots,M\}$ satisfy (\ref{assumption_new}) for all $r \geq r_0$. By the triangle inequality in $L^2(B_r^+)$ and the Poincar\'{e} inequality on $B_r^+$ with homogeneous Dirichlet boundary conditions on $\bhs\cap B_r$, writing $\phi_d - \tilde{\varphi}^{\Sigma}_{M}=(1-\sum_{m=-1}^M\chi_m)\phi_d -\varphi^{\Sigma}_{M}$ we get
\begin{align}
\nonumber
&\frac{1}{r} \left(\fint_{B_{r}^+} \left| \phi_d - \tilde{\varphi}^{\Sigma}_{M}\right|^2
+\left|\sigma_d-\psi_{,M} \right|^2 \,dx\right)^{1/2}
\\&
\label{sublinearinduct1}
\leq \frac{1}{r}\left(\fint_{B_r^+}|(\phi_d, \sigma_d)|^2 \,dx \right)^{1/2}
+\frac{1}{r} \left(\fint_{B_r^+} |\psi_{,M}|^2 \,dx \right)^{1/2} 
\\&~~~\nonumber
+C_P \left(\fint_{B_r^+} | \nabla \varphi^{\Sigma}_{M}|^2 \,dx \right)^{1/2}.
\end{align}
Notice that for $r \geq r_0$ Lemma \ref{2.3} yields
\begin{align}
\nonumber
&\frac{1}{r} \left( \fint_{B_r^+} |\psi_{,M}|^2     \right)^{1/2}
\\&
\leq 2 d^2  C_3 \sum_{n=-1}^{\infty} \min\bigg\{ 1, \frac{r_0 2^{n+1}}{r} \bigg\} \fir{r_0 2^{n+1}}{\nabla \tilde{\varphi}^{\Sigma}_{M}}.
\label{sublinearinduct2}
\end{align}
Using our assumption that the $\varphi_m$ with $m\in \{-1,\ldots,M\}$ satisfy (\ref{assumption_new}) for any $r\geq r_0$ -- and therefore in particular for $r:=r_0 2^{n+1}$ -- gives that
\begin{align}
\nonumber
&\sum_{n=-1}^{\infty} \min\bigg\{1,\frac{r_0 2^{n+1}}{r}\bigg\} \fir{r_0 2^{n+1}}{\nabla \varphi^{\Sigma}_M}
\\&
\nonumber
\leq C_2 \sum_{m=-1}^M \sum_{n=-1}^{\infty}
\min\big\{1, 2^{\frac{d(m-n)}{2}}\big\} \delta_{r_0 2^{m+1}}^{1/3}
\\&
\label{piece1}
\leq C_2 \sum_{m=-1}^M \Big(m + 1 + \frac{1}{1-2^{-d/2}}\Big) \delta_{r_0 2^{m+1}}^{1/3}. 
\end{align}
Furthermore, we may use that $\chi_m$ is supported in $B_{r_0 2^{m+1}}^+ \setminus B_{r_0 2^{m-1}}^+$ for $m\neq -1$ and \eqref{chiphiformula} (applied with $r:=r_0 2^{n+1}$) to get that
\begin{align}
\nonumber
&\sum_{m=-1}^M \displaystyle\sum_{n=-1}^{\infty} \min\bigg\{1, \frac{ r_0 2^{n+1}}{r} \bigg\} \fir{r_0 2^{n+1}}{\nabla (\chi_m \phi_d)}
\\& \nonumber
\leq \sum_{m=-1}^M \sum_{n=m-1}^{\infty} \fir{r_02^{n+1}}{\nabla(\chi_m \phi_d)}
\\&\nonumber
\leq C_1 \sum_{m=-1}^M \sum_{n=m-1}^{\infty} 2^{\frac{d(m-n)}{2}}\delta_{r_0 2^{m+1}}^{1/3}
\\&
\leq C_1 \sum_{m=-1}^M \frac{2^{d/2}}{1-2^{-d/2}} \delta_{r_0 2^{m+1}}^{1/3}.
\label{piece2}
\end{align}
Then, continuing (\ref{sublinearinduct2}) with (\ref{piece1}) and (\ref{piece2}) yields
\begin{align*}
&\frac{1}{r} \left(\fint_{B_r^+} |\psi_{,M}|^2 \,dx \right)^{1/2}
\\&
\leq 2 d^2  C_3 (C_1 + C_2) \sum_{m=-1}^M \bigg(m + 1 + \frac{2^{d/2}}{1-2^{-d/2}} \bigg) \delta_{r_0 2^{m+1}}^{1/3}
\\&
\leq 2 d^2  C_3 (C_1 + C_2) \sum_{k=m_0}^{M+m_0+1} \bigg(k + \frac{2^{d/2}}{1-2^{-d/2}}\bigg) \delta_{2^k}^{1/3}.
\end{align*}
To treat the other term of (\ref{sublinearinduct1}) we again use (\ref{assumption_new}), which gives
\begin{align*}
&\left(\fint_{B_r^+} |\nabla \varphi^{\Sigma}_{M} |^2 \,dx \right)^{1/2}
\\
&\leq C_2 \sum_{m=-1}^{M}\min \bigg\{1, \bigg(\frac{r_0 2^{m+1}}{r}\bigg)^{d/2}\bigg\}\delta_{r_0 2^{m+1}}^{1/3}
\\&
\leq C_2 \sum_{k=m_0}^{M+m_0+1} \delta_{2^k}^{1/3}.
\end{align*}
So, for $ r \geq r_0$ we arrive at
\begin{align}
\nonumber
&\frac{1}{r} \left(\fint_{B_{r}^+} \left| \phi_d -  \tilde{\varphi}^{\Sigma}_{M}
\right|^2
+\left|\sigma_d-\psi_{,M} \right|^2 \vphantom{\fint_{B_{r}^+}} \,dx \right)^{1/2} 
\\& \nonumber
\leq 2\delta_r + 2 d^2 C_3 (C_1+C_2) \sum_{k=m_0}^{\infty} \bigg(k + \frac{2^{d/2}}{1-2^{-d/2}} \bigg) \delta_{2^{k}}^{1/3}
\\&~~~
+ C_P C_2 \sum_{k=m_0}^{\infty} \delta_{2^{k}}^{1/3}.
\label{sublinearinduct5}
\end{align}
As a consequence of the estimate \eqref{sublinearinduct5}, our assumption (\ref{condition}) allows us to choose $r_0=2^{m_0}$ large enough -- independently of $M$ -- such that for $(\tp_{d,M}, \ts_{d,M})$ the estimate (\ref{smallness}) is satisfied for $\alpha=1/2$ and $r\geq r_0$.

Thus, we infer the estimate \eqref{assumption_new} for $\varphi_{M+1}$: The case $\frac{r_0 2^{(M+1)+1}}{r}\leq 16$ has already been treated in Lemma \ref{first}; it just remains to extend the estimate to the case $\frac{r_0 2^{(M+1)+1}}{r}>16$. As $(\tp_{d,M}, \ts_{d,M})$ is a half-space-adapted corrector on $B_R^+$ with $R:=r_0 2^M$ which satisfies \eqref{smallness} for $\alpha=1/2$ and $r\geq r_0$, Theorem \ref{excessbound} is applicable and yields the mean-value property \eqref{mvp} for $a$-harmonic functions on $B_{r_0 2^M}^+$ with homogeneous Dirichlet boundary conditions on $\bhs\cap B_{r_0 2^M}$. Since $\varphi_{M+1}$ is indeed $a$-harmonic in $B_{r_0 2^M}^+$ with homogeneous Dirichlet boundary conditions on $\bhs \cap B_{r_0 2^M}$, we deduce for $r\in [r_0,r_0 2^M]$ using in the second step the estimate \eqref{energy1formula} for $r:=r_0 2^{M}$
\begin{align*}
&\fir{r}{\nabla \varphi_{M+1} } \leq C_{Mean} \fir{r_0 2^{M}}{\nabla \varphi_{M+1} }
\\&
\leq  C_{Mean} C_1 2^d \delta_{r_0 2^{(M+1)+1}}^{1/3}.
\end{align*}
This shows \eqref{assumption_new} for $\varphi_{M+1}$ and $r\in [r_0,r_0 2^M]$.
\end{proof}

\subsection{Step 4 -- Passage to the limit $M\rightarrow\infty$}

\begin{proof}[Proof of Theorem \ref{existenceofcorrectors}]
Let the $L_m$ be chosen as in Lemma \ref{first}. Let $r_0=2^{m_0}$ be chosen as in Lemma \ref{induct}. By Lemma \ref{induct}, the estimate \eqref{assumption_new} then holds for all $m\geq -1$ (the start of the induction -- i.\,e.\ \eqref{assumption_new} for $m=-1$ -- is provided by Lemma \ref{first}).

For $i\neq d$ we then choose $\phi_i^\h:=\phi_i|_{\hs}$ and $\sigma_{ijk}^\h:=\sigma_{ijk}|_{\hs}$. By our assumption \eqref{condition}, we therefore have to verify the assertion on sublinear growth iv) in our theorem only for $\phi_d^\h$ and $\sigma_{d}^\h$.
{~\\ \bf Part 1}: The corrector $\phi_d^\h$.\nopagebreak\\
We first show that the series $\sum_{m=-1}^\infty \varphi_m$ converges absolutely in $H^1(B_r^+)$ for all $r \geq r_0$. By the Poincar\'{e} inequality for functions in $H^1(B_r^+)$ with homogeneous Dirichlet boundary conditions on $\bhs\cap B_r$, it suffices to calculate (using (\ref{assumption_new}))
\begin{align*}
&\sum_{m=-1}^{\infty} \fir{r}{\nabla \varphi_{m}}
\leq C_2 \sum_{m=-1}^{\infty} \delta_{r_0 2^{m+1}}^{1/3}
\leq C_2 \sum_{k=m_0}^{\infty} \delta_{2^k}^{1/3}
\end{align*}
and to use the summability of the $\{\delta_{2^{k}}^{1/3}\}_{k}$ (see \eqref{condition}).

Again, combining (\ref{assumption_new}) with the Poincar\'{e} inequality yields for the sum $\varphi:=\sum_{m=-1}^\infty \varphi_m =\lim_{M\rightarrow \infty} \varphi_M^\Sigma$
\begin{align}
\nonumber
\frac{1}{r} \fir{r}{ \varphi } 
&\leq
\sup_M \frac{1}{r} \fir{r}{ \varphi_M^{\Sigma}}
\\&
\label{final1}
\leq C_PC_2 \sum_{k=m_0}^{\infty} 
\min \bigg\{1, \left(\frac{2^{k}}{r}\right)^{d/2}\bigg\}
\delta_{2^{k}}^{1/3}
\end{align}
for all $r \geq r_0$.

Next, we show that $\{ \sum_{m=-1}^{M} \chi_m  \phi_d \}_M$ forms a Cauchy sequence in $H^1(B_r^+)$ for all $r \geq r_0$. Using the fact that $\chi_m \phi_d$ vanishes outside of $B_{r_0 2^{m+1}}\setminus B_{r_0 2^{m-1}}$ (except for $m=-1$, for which $\chi_{-1}\phi_d$ vanishes outside of $B_{r_0}$), the Poincar\'{e} inequality for functions in $H^1(B_{r_02^{m+1}})$ that vanish on $\partial B_{r_02^{m+1}}\cap \hs$ yields that for any $r>0$ 
\begin{align}
\label{beforenabla}
\left(\int_{B_r} |\chi_{m} \phi_d|^2 \,dx \right)^{1/2} \lesssim r \left(\int_{B_r}|\nabla (\chi_{m}\phi_d)|^2 \,dx \right)^{1/2}.
\end{align}
Using \eqref{chiphiformula} and again $\operatorname{supp}\chi_m \subset B_{r_0 2^{m+1}}\setminus B_{r_0 2^{m-1}}$, we see that
\begin{align}
\nonumber
\sum_{m=-1}^{\infty} \left(\fint_{B_r} |\nabla (\chi_m \phi_d)|^2 \,dx \right)^{1/2} &
\leq 2^d C_1 \sum_{m=-1}^{\infty} 
\min \bigg\{1, \bigg(\frac{r_0 2^{m+1}}{r}\bigg)^{d/2}\bigg\}\delta_{r_0 2^{m+1}}^{1/3}
\\
\label{nabla}
&\leq 2^d C_1 \sum_{k=m_0}^{\infty} 
\min \bigg\{1, \bigg(\frac{2^{k}}{r}\bigg)^{d/2}\bigg\}
\delta_{2^{k}}^{1/3}.
\end{align}
So,  $\left\{ \sum_{m=-1}^{M} \chi_m  \phi_d \right\}_M$ forms a Cauchy sequence in $H^1(B_r^+)$.

The function $\tilde{\varphi} := \varphi + \sum_{m=-1}^{\infty} \chi_m \phi_d = \lim_{M\rightarrow \infty} \tilde \varphi_M^\Sigma$ is a weak solution of the problem (\ref{TildeVarphiProblem}): \eqref{TildeVarphiProblem} is satisfied on $B_r$ by all $\tilde \varphi_M^\Sigma$ for which $r_0 2^{M} \geq r$ holds. Thus, \eqref{TildeVarphiProblem} carries over to the limit $M\rightarrow \infty$ for arbitrarily big radii $r$.  Therefore \eqref{TildeVarphiProblem} holds globally for the limit $\tilde\varphi$, which entails that $\phi_d^\h=\phi_d-\tilde\varphi$ solves \eqref{halfspacecorrector}.

By (\ref{final1}), \eqref{beforenabla}, and \eqref{nabla}, our assumption (\ref{condition}) implies that $\tilde{\varphi}$  and, therefore,  $\phi_d^{\h} = \phi_d - \tilde{\varphi}$ are sublinear in the sense
\begin{align*}
\lim_{r\rightarrow \infty} \frac{1}{r} \left(\fint_{B_r} |\phi_d^\h|^2 \,dx\right)^{1/2} =0.
\end{align*}

\smallskip
\noindent{\bf Part 2:} The vector potential $\sigma_d^\h$.\\
We now show that $\{\psi_{jk,M}\}_{M}$ forms a Cauchy sequence in $L^2(B_r^+)$ for all $ r \geq r_0$; furthermore, we show that the limit $\psi_{jk}$ has sublinear growth.
To this aim, observe that the differences $v_{j,M+1}^n - v_{j,M}^n$ are weak solutions to the problem
\begin{subequations}
\label{energy3}
\begin{align}
-\Delta (v_{j,M+1}^n - v_{j,M}^n )
&=
-\di ( \eta_n x_j a  \nabla (\varphi_{M+1} + \chi_{M+1} \phi_d ))
&&\textrm{in } \hs,
\\
v_{j,M+1}^n - v_{j,M}^n &= 0 &&\textrm{if } j \neq d \textrm{ on } \bhs,
\\
\partial_d (v_{d,M+1}^n - v_{d,M}^n )&= 0  &&\textrm{on } \bhs.
\end{align}
\end{subequations}
To shorten the subsequent computations, let us use the convention $v_{j,-2}^n\equiv 0$ and $b_{j,-2}^n=0$; then \eqref{energy3} holds also for $M=-2$.

Estimating analogously to the proof of Lemma \ref{2.2} -- note that the only difference between the equation satisfied by $v_{j,M}^n$ and the equation satisfied by $v_{j,M+1}^n-v_{j,M}^n$ is the right-hand side -- , we deduce that for any $r\geq r_0$
\begin{align*}
&\frac{1}{r}\fir{r}{\partial_k \left(v_{j,M+1}^n - v_{j,M}^n - (b_{j,M+1}^n - b_{j,M}^n)\cdot x  \right)}
\\&\nonumber
\leq C_3 \min\bigg\{1,\frac{r_0 2^{n+1}}{r}\bigg\} \fir{r_0 2^{n+1}}{\nabla (\varphi_{M+1} + \chi_{M+1} \phi_d)}.
\end{align*}
Taking the sum with respect to $n$, we deduce that the limits $v_{j,M}$ of the series $\sum_{n=-1}^\infty (v_{j,M}^n - b_{j,M}^n \cdot x)$ satisfy
\begin{align*}
&\frac{1}{r}\fir{r}{\partial_k (v_{j,M+1} - v_{j,M})}
\\&\nonumber
\leq C_3 \sum_{n=-1}^\infty \min\bigg\{1,\frac{r_0 2^{n+1}}{r}\bigg\} \fir{r_0 2^{n+1}}{\nabla (\varphi_{M+1} + \chi_{M+1} \phi_d)}.
\end{align*}
Taking the sum with respect to $M$ and estimating the right-hand side by the inequality (\ref{assumption_new}) and the estimate \eqref{chiphiformula} -- both inequalities applied with $r$ replaced by $r_0 2^{n+1}$ and $m$ replaced by $M+1$ -- (note again that $\chi_{M+1}\phi_d$ vanishes on $B_{r_0 2^{n+1}}^+$ in case $r_0 2^{M+1-1}\geq r_0 2^{n+1}$ and $M+1\neq -1$), we infer
\begin{align}
\label{QuantitativeSublinear}
&\frac{1}{r}\sum_{M=-2}^\infty \fir{r}{\partial_k (v_{j,M+1} - v_{j,M})}
\\&
\nonumber
\leq 
C_3(C_2+ 2^d C_1) \sum_{n=-1}^\infty \min\bigg\{ 1, \frac{r_0 2^{n+1}}{r}\bigg\}  \sum_{M=-2}^{\infty} 
\min \big\{1,2^{d(M+1-n)/2}\big\}\delta_{r_0 2^{M+1+1}}^{1/3}.
\end{align}
Now, by this estimate, it is sufficient to show
\begin{align}
\label{SumConverges}
\sum_{n=-1}^\infty \sum_{M=-2}^{\infty} \min \big\{1,2^{d(M+1-n)/2}\big\}\delta_{r_0 2^{M+1+1}}^{1/3}<\infty
\end{align}
in order to obtain both the Cauchy sequence property of $\nabla v_{j,M}$ in $L^2(B_r)$ and the sublinearity property
\begin{align*}
\lim_{r\rightarrow \infty} \frac{1}{r} \fir{r}{\partial_k v_{j}}
\leq
\lim_{r\rightarrow \infty} \sup_{M\geq -1} \frac{1}{r} \fir{r}{\partial_k v_{j,M}}
=0.
\end{align*}
Note that by $\psi_{jk}:=\partial_k v_j-\partial_j v_k$ and $\sigma_{djk}^\h=\sigma_{djk}-\psi_{jk}$, this estimate then directly implies the desired result
\begin{align*}
\lim_{r\rightarrow \infty} \frac{1}{r} \fir{r}{\sigma_d^\h} =0.
\end{align*}
Furthermore, the $\psi_{jk,M}$ are solutions to the equation \eqref{equationforpsi}. Since we can pass to the limit $M\rightarrow \infty$ in the weak formulation of \eqref{equationforpsi} for any smooth compactly supported test function, this shows that the limit $\sigma_{djk}^\h:=\lim_{M\rightarrow \infty} (\sigma_{djk}-\psi_{jk,M})$ solves the equation \eqref{halfspacepotential}.

However, to see that \eqref{SumConverges} holds, we just need to estimate
\begin{align*}
&\sum_{n=-1}^\infty \sum_{M=-2}^{\infty} \min \big\{1,2^{d(M+1-n)/2}\big\}\delta_{r_0 2^{M+1+1}}^{1/3}
\\&
\leq
\sum_{M=-2}^{\infty} \Big(M+2+\frac{1}{1-2^{-d/2}}\Big)
\delta_{r_0 2^{M+1+1}}^{1/3}
=
\sum_{k=m_0}^{\infty} \Big(k-m_0+\frac{1}{1-2^{-d/2}}\Big)
\delta_{2^k}^{1/3}
\end{align*}
and use the summability property \eqref{condition}. This finishes the proof of our theorem.
\end{proof}

\section{Proofs of the Regularity Results Theorem \ref{excessbound} and Corollary \ref{Liouville}}

In the proof of Theorem \ref{excessbound}, we shall need the following Caccioppoli inequality.
\begin{customlemma}{4.1}
\label{caclemma}
Let $a$ be a coefficient field satisfying the ellipticity and boundedness assumptions \eqref{Elliptic}. For any $a$-harmonic function $u$ on $B_R^+$ subject to homogeneous Dirichlet boundary conditions on $\bhs \cap \partial B_R^+$, the estimate
\begin{align}
\label{cac}
\fir{R/2}{\nabla u} \lesssim \frac{1}{R} \fir{R}{u}.
\end{align}
holds.
\end{customlemma}
\begin{proof}
Testing the equation
\begin{align*}
 -\di (a \nabla u) =0 \quad \textrm{in $B_R^+$}
\end{align*}
with $\eta^2 u$, where $\eta$ is a radial cut-off with $\eta\equiv 1$ in $B_{R/2}$, $\eta \equiv 0$ outside of $B_R$, $0\leq \eta\leq 1$ everywhere, and $|\nabla\eta|\leq \frac{3}{R}$, we get
\begin{align*}
\ipts{B_R^+}{\eta^2 \nabla u \cdot a \nabla u + 2 \eta u \nabla \eta \cdot a \nabla u}
=0.
\end{align*}
Note that the boundary terms vanish as $\eta^2 u$ is zero on $\partial B_R^+$. Using the uniform ellipticity of $a$ and Young's inequality allows us to write
\begin{align*}
 \lambda \ipts{B_R^+}{\eta^2 |\nabla u|^2} \leq  2 \ipts{B_R^+}{ |\eta u \nabla \eta \cdot  a\nabla u|} \leq   \ipts{B_R^+}{\frac{\lambda}{2} \eta^2 |\nabla u|^2 + \frac{2}{\lambda}|\nabla \eta|^2  u^2}.
\end{align*}
The properties of $\eta$ finish the argument.
\end{proof}
The following classical regularity properties of constant-coefficient elliptic equations will play a crucial rule in the derivation of the excess-decay estimate.
\begin{customlemma}{4.2}
Let $v$ be a weak solution to the constant-coefficient equation $-\nabla \cdot (a_{hom}\nabla v)=0$ in $B_{\Rp}^+$ with homogeneous Dirichlet boundary conditions on $\partial B_{\Rp}^+ \cap \bhs$, where $a_{hom}$ is a positive definite matrix. Then there exists some $\beta=\beta(d,\lambda)>0$ such that for any positive $\rho\leq \frac{1}{2}R'$ and any positive $r\leq \frac{1}{2}R'$ the following estimates hold:
\begin{subequations}
\label{innerregularity}
\begin{align}
\label{innerregularitya}
r^2 \displaystyle\sup_{B_r^+}|\nabla^2 v|^2 & \lesssim \left(\displaystyle\frac{r}{\Rp}\right)^2 \firs{\Rp}{\nabla v},
\\
\label{innerregularityb}
\irt{B_{\Rp}^+ \setminus B_{\Rp-2\rho}^+}{\nabla v} & \lesssim R' \left(\frac{\rho}{R'} \right)^{\beta}\int_{\partial B_{\Rp}^+}|\nabla^{tan} v |^2 \,dS,
\\
\label{innerregularityc}
\displaystyle\sup_{B_{\Rp-\rho}^+}(|\nabla^2 v|^2+\frac{1}{\rho^2}|\nabla v|^2)& \lesssim \displaystyle\frac{1}{\rho^2}\left( \displaystyle\frac{\Rp}{\rho}\right)^d \firs{\Rp}{\nabla v}.
\end{align}
\end{subequations}
\end{customlemma}
\begin{proof}
For the third estimate, notice that if $x^{\prime} \in S$ where $S= B_{\Rp-\rho}^+ \cap \{x_d\geq \frac{\rho}{2}\}$ then $v$ is $a_{hom}$-harmonic on
$B_{\rho/2}(x^{\prime})$. Therefore, for these $x^{\prime}$ we have the inner regularity estimate
\begin{align}
\label{innerregularitywb}
\displaystyle\sup_{y \in B_{\rho/4}(x^{\prime})} \rho^2|\nabla^2 v(y)|^2 + \displaystyle\sup_{y \in B_{\rho/4}(x^{\prime})}  |\nabla v(y)|^2
\lesssim \displaystyle\frac{1}{\rho^d} \displaystyle\int_{B_{\rho/2}(x^{\prime})}|\nabla v|^2\,dx,
\end{align}
which follows by an iterated use of the Caccioppoli inequality on balls to derive an $H^k$ estimate for $k$ large enough and a subsequent use of the Sobolev embedding.

For $x^{\prime} \in S_{\partial}$, where $S_{\partial} = \bhs \cap B_{\Rp-\rho}^+$, we get an analogous estimate for half-balls: In this case, the result can also be shown by proving $H^k$ regularity estimates for $k$ large enough followed by the Sobolev embedding. The derivation of $H^k$-type regularity estimates is again standard: One may proceed by repeatedly using the Caccioppoli estimate for $v$ and its tangential (higher) derivatives $\partial_{i_1} \ldots \partial_{i_{k-1}} v$ with $i_1,\ldots,i_{k-1}\neq d$. To obtain estimates on higher derivatives which involve multiple derivatives in the normal direction $e_d$ -- only estimates for derivatives containing a single normal derivative are provided by the aforementioned applications of the Caccioppoli inequality -- one directly uses the equation satisfied by $v$.
Thus, for $x^{\prime} \in S_{\partial}$ we have
\begin{align}
\label{innerregularityhb}
\displaystyle\sup_{y \in B_{\rho/2}^+(x^{\prime})} \rho^2|\nabla^2 v(y)|^2 + \displaystyle\sup_{y \in B_{\rho/2}^+(x^{\prime})}  |\nabla v(y)|^2
 \lesssim \displaystyle\frac{1}{\rho^d} \displaystyle\int_{B_{\rho}^+(x^{\prime})}|\nabla v|^2 \,dx.
\end{align}
The estimate (\ref{innerregularitya}) is an immediate consequence of (\ref{innerregularityhb}) with $\rho:=R'$ and $x'=0$. To obtain (\ref{innerregularityc}) let 
\begin{align*}
s &= \displaystyle\sup_{x^{\prime}\in S} \displaystyle\sup_{y \in B_{\rho/4}(x^{\prime})} (|\nabla^2 v(y)|^2+\frac{1}{\rho^2}|\nabla v(y)|^2),
\\
s_{\partial} &= \displaystyle\sup_{x^{\prime}\in S_{\partial}} \displaystyle\sup_{y \in B^+_{\rho/2}(x^{\prime})} (|\nabla^2 v(y)|^2+\frac{1}{\rho^2}|\nabla v(y)|^2).
\end{align*}
Using (\ref{innerregularitywb}) and (\ref{innerregularityhb}), we may then write
\begin{align*}
&\displaystyle\sup_{x \in B_{\Rp-\rho}^+}(|\nabla^2 v|^2+\frac{1}{\rho^2}|\nabla v|^2) 
\leq \displaystyle\max \{s, s_{\partial}\}
\\
&\lesssim \displaystyle\sup_{x^{\prime} \in S \cup S_{\partial}}\displaystyle\frac{1}{\rho^{d+2}}  \displaystyle\int_{B_{\rho}(x^{\prime})\cap \hs}|\nabla v|^2\,dx
\lesssim \displaystyle\frac{1}{\rho^2} \left(\displaystyle\frac{\Rp}{\rho} \right)^d  \firs{\Rp}{\nabla v},
\end{align*}
finishing the proof of \eqref{innerregularityc}.

Finally, for the inequality \eqref{innerregularityb} we first extend $v$ to $B_{R'}$ by odd-reflection. The extended $v$ satisfies the elliptic equation
\begin{align*}
~~~~~~~~~-\nabla \cdot (\tilde a_{hom} \nabla v)=0 && \text{in }B_{R'}
\end{align*}
with
\begin{align*}
(\tilde a_{hom})_{ij}=
\begin{cases}
(a_{hom})_{ij} &\text{ for }x_d>0,
\\
(a_{hom})_{ij}&\text{ for }x_d<0 \text{ and }i\neq d, j\neq d,
\\
-(a_{hom})_{ij}&\text{ for }x_d<0 \text{ and }i=d, j\neq d,
\\
-(a_{hom})_{ij}&\text{ for }x_d<0 \text{ and }i\neq d, j= d,
\\
(a_{hom})_{ij}&\text{ for }x_d<0 \text{ and }i=j=d.
\end{cases}
\end{align*}
If we then let $\bar{v}$ be the harmonic extension of $v|_{\partial B_{R'}}$ to $B_{\Rp}$, we have the estimate $||\nabla \bar{v}||_{L^{2/(1-\beta)}(B_{\Rp})} \lesssim  {R'}^{1/2-d\beta/2} ||\nabla^{tan}v||_{L^2(\partial B_{\Rp})}$, provided that $\beta>0$ is not too large. Furthermore, Meyers' estimate \cite{Meyers} states that for any $\beta>0$ small enough (depending on $d$ and $\lambda$), the solution $v-\bar v$ to the equation
\begin{align*}
-\nabla \cdot(\tilde a_{hom} \nabla (v-\bar v))
&= \nabla \cdot (\tilde a_{hom} \nabla \bar v)
&&\text{in }B_{\Rp},
\\
v-\bar v &=0&&\text{on }\partial B_{\Rp}
\end{align*}
satisfies the bound $||\nabla (v-\bar v)||_{L^{2/(1-\beta)}(B_{\Rp})}\lesssim ||\tilde a_{hom}\nabla \bar v||_{L^{2/(1-\beta)}(B_{\Rp})}$.
%$v \in H^1(B_{\Rp})$ and $v-\bar{v} \in H_0^1(B_{\Rp})$ then $||\nabla v||_{L^p(B_{\Rp})} \lesssim||\nabla \bar{v}||_{L^p(B_{\Rp})}$ for $p = \frac{2d}{d-1}$ \cite{GilbargTrudinger}. 
Combining this estimate with the bound on $\bar v$ yields that 
\begin{align*}
||\nabla v ||_{L^{2/(1-\beta)}(B_{\Rp}^+)} \lesssim {R'}^{1/2-d\beta/2} ||\nabla^{tan} v||_{L^2(\partial B_{\Rp}^+)}.
\end{align*}
%Here, we have used that $\nabla^{tan}v=0$ in $\partial \hs$.
It then follows by H\"older's inequality that
\begin{align*}
&\left(\int_{B_{\Rp}^+ \setminus B_{\Rp-2\rho}^+} |\nabla v|^2 \,dx \right)^{1/2} \leq |B_{\Rp}^+ \setminus B_{\Rp-2\rho}^+|^{\beta/2} \left( \int_{B_{\Rp}^+}|\nabla v|^{2/(1-\beta)} \,dx\right)^{(1-\beta)/2}
\\&
\lesssim (R^{\prime})^{1/2-\beta/2} \rho^{\beta/2} \left( \int_{\partial B_{\Rp}^+}|\nabla^{tan} v |^2 \,dS \right)^{1/2},
\end{align*}
concluding the proof of \eqref{innerregularityb}.
\end{proof}

We now turn to the proof of the excess-decay estimate.
\begin{proof}[Proof of Theorem \ref{excessbound}]
~\\
\noindent\textbf{Step 1}:\\
In the first step of the proof, we show that for each $r<R$ there exists $b \in \mathbb{R}$ such that the estimate
\begin{align}
\nonumber
&\firs{r}{\nabla u - b (e_d +\nabla \tp_d)}
\\
&\lesssim  \bigg( \left(\displaystyle\frac{r}{R}\right)^2 \left(1 + \delta^2 \right)  
+\left( \displaystyle\frac{R}{r} \right)^d \delta^{2\beta/(d+2+\beta)} \vphantom{\left(\displaystyle\frac{r}{R}\right)^2} \bigg)
\firs{R}{\nabla u}
\label{iterate}
\end{align}
is valid, with the abbreviation
\begin{align*}
\delta:=\max \big\{ \delta^{\h}_{2r}, \delta^{\h}_{R}\big\}.
\end{align*}
In the proof, for convenience we make use of the \emph{Einstein summation convention}, i.\,e.\ whenever an index appears twice in an expression, summation with respect to the index is implied.

Note that for $r \in \left[\frac{R}{4}, R \right]$ the estimate trivially holds for $b = 0$. It is therefore sufficient to show (\ref{iterate}) for $r\leq R/4$. To do this, we first choose a radius $\Rp \in ( \frac{R}{2}, R )$ such that 
\begin{equation}\label{choice}
 \irt{\partial B_{\Rp}^+}{\nabla^{tan} u} \lesssim \displaystyle\frac{1}{R}\displaystyle\int_{B_{R}^+ \setminus B_{R/2}^+} | \nabla u|^2 \,dx \lesssim \displaystyle\frac{1}{R}\irs{R}{\nabla u}.
\end{equation}
We know that such a radius exists by writing the middle integral in polar coordinates and using that $\nabla^{tan} u=0$ on $\bhs\cap B_R$.

Let $v$ be the $a_{hom}$-harmonic function that coincides with $u$ on $\partial B_{R^{\prime}}^+$. To show the estimate (\ref{iterate}) we compare $\nabla u$ to $\nabla v$ corrected as suggested by the two-scale expansion \eqref{TwoScaleExpansion}. Notice that, due to the boundary conditions of $v$, we know that $\nabla v(0)$ only has a normal component. This observation allows us to write
\begin{align}
\nonumber
&\irt{B_r^+}{\nabla u - \partial_d v(0) (e_d +\nabla \tp_d)}
\\
\label{Taylor}
&\lesssim \irt{B_r^+} {(\nabla v - \nabla v(0))(\textrm{id} + \nabla \tp)}
+\irt{B_r^+}{\nabla u - \partial_i v (e_i + \nabla \tp_i)}.
\end{align}
Notice that the second term on the right hand side is the gradient of the ``homogenization error'' coming from the ansatz for $v$ given by the two-scale expansion; see \eqref{TwoScaleExpansion}.
To estimate this term, we first derive an estimate for
\begin{align*}
w := u - (v + \eta \tp_i\partial_i v),
\end{align*}
where $\eta$ is a cut-off with $0\leq \eta\leq 1$, $\eta\equiv 1$ in $B_{\Rp-2\rho}^+$, $\eta\equiv 0$ outside of $B_{\Rp-\rho}^+$, and $|\nabla \eta|\leq \frac{C}{\rho}$.
We will later optimize the width of the boundary-layer introduced by $\rho$, but for the moment we only assume that $0<\rho \leq \frac{1}{4} \Rp$. The function $w$ satisfies the equation
\begin{align}
\label{equnforw}
 -\di (a \nabla w) = \di ((1-\eta) (a-a_{hom})\nabla v +(\tp_i a -\ts_i)\nabla(\eta \partial_i v)) \quad \text{in }B_{\Rp}^+.
\end{align}
To see this, one uses that $u$ is $a$-harmonic, that $\tp_i$ solves the corrector equation \eqref{CorrectorPhiEquation} on $B_{\Rp}^+$, and the defining property \eqref{def_sigma} of $\ts$, which gives
\begin{align*}
&-\di (a \nabla w)
\\&
=\di \left( a \nabla v + \eta \partial_i v a  \nabla \tp_i \right) + \di (\tp_i a  \nabla (\eta \partial_i v))
\\&
= \di \left((1-\eta)a\nabla v + \eta \partial_i v a(e_i+\nabla \tp_i) \right)+ \di (\tp_i a  \nabla (\eta \partial_i v))
\\&
= \di ((1-\eta)a\nabla v) + \nabla(\eta \partial_i v)\cdot a(e_i+\nabla \tp_i) + \di (\tp_i a  \nabla (\eta \partial_i v))
\\&
= \di ((1-\eta)(a-a_{hom}) \nabla v) +\nabla(\eta \partial_i v)\cdot(a(e_i+\nabla \tp_i)-a_{hom}e_i)
\\&~~~
+ \di (\tp_i a  \nabla (\eta \partial_i v))
\\&
= \di ((1-\eta)(a-a_{hom}) \nabla v) +\nabla(\eta \partial_i v)\cdot ( \di \ts_i )
+ \di (\tp_i a  \nabla (\eta \partial_i v)).
\end{align*}
To complete the calculation, we use the skew-symmetry of the vector potential $\ts_{ijk}$ in the form
$\nabla(\eta \partial_i v)\cdot ( \di \ts_i ) = -\nabla \cdot (\ts_{i} \nabla (  \eta \partial_i v))$.

Notice that, due to the cut-off $\eta$, the boundary conditions of $\tp_d$, and the boundary conditions of $v$, $w$ satisfies homogeneous Dirichlet boundary conditions on $\partial B_{R^{\prime}}^+$. Therefore, the standard energy estimate for the equation (\ref{equnforw}) reads
\begin{align*}
&\ir{\Rp}{\nabla w}
\\&
\leq \frac{1}{\lambda} \ir{\Rp}{(1-\eta)(a-a_{hom})\nabla v + (\tp_i a - \ts_i)\nabla(\eta \partial_i v)}.
\end{align*}
The boundedness of $a$ and $a_{hom}$ and the properties of $\eta$ then imply
\begin{align}
\nonumber
&\irs{\Rp-2\rho}{\nabla u - \partial_i v(e_i + \nabla \tp_i)}
\\&
\lesssim \irt{B_{\Rp}^+\setminus B_{\Rp-2\rho}^+}{\nabla v}
+\int_{B_{\Rp-\rho}^+} |(\tp,\ts)|^2 (|\nabla^2 v|^2+\frac{1}{\rho^2}|\nabla v|^2)\,dx.
 \label{energyestimate}
\end{align}
Due to the conditions that we have placed on $r$, $\rho$, and $R^{\prime}$ we have $r \leq \Rp-2\rho$. Therefore the second term on the right hand side of (\ref{Taylor}) can be estimated by the formula (\ref{energyestimate}).  This yields
\begin{align}
\nonumber
&\irt{B_r^+}{\nabla u - \partial_d v(0) (e_d +\nabla \tp_d)}
\\&\nonumber
\lesssim \displaystyle\int_{B_r^+}{|\nabla v-\nabla v (0)|^2}|\textrm{id}+\nabla \tp|^2 \,dx
\\&~~~\nonumber
+\irt{B_{\Rp}^+-B_{\Rp-2\rho}^+}{\nabla v} +\displaystyle\int_{B_{\Rp-\rho}^+} |(\tp, \ts)|^2 (|\nabla^2 v|^2+\frac{1}{\rho^2}|\nabla v|^2)\,dx
\\&\nonumber
\leq r^2 \displaystyle\sup_{B_r^+}|\nabla^2 v|^2 \displaystyle\int_{B_r^+}{|\textrm{id}+\nabla \tp|^2}\,dx
\\&~~~
+\irt{B_{\Rp}^+\setminus B_{\Rp-2\rho}^+}{\nabla v} +\displaystyle\sup_{B_{\Rp-\rho}^+}(|\nabla^2 v|^2+\frac{1}{\rho^2}|\nabla v|^2) \displaystyle\int_{B_R^+} |(\tp, \ts)|^2 \,dx.
\label{step3}
\end{align}
To further process this estimate, we exploit that $v$ solves the constant-coefficient equation $-\nabla \cdot (a_{hom}\nabla v)=0$ in $B_{\Rp}^+$ with homogeneous Dirichlet boundary conditions on $\bhs\cap B_{\Rp}$; thus the estimates \eqref{innerregularity} are available. Furthermore, notice that the difference $v-u$ solves
\begin{align*}
~~~~~~~~~~-\di(a_{hom} \nabla (v-u)) &= \di(a_{hom}\nabla u) &&\textrm{in }\quad B_{\Rp}^+,
\\
v-u &= 0 &&\textrm{on }\quad \partial B_{\Rp}^+.
\end{align*}
Testing this equation with $v-u$ and using Young's inequality yields
\begin{equation}
\label{apriori}
\irs{\Rp}{\nabla v}
\leq 2\irs{\Rp}{\nabla u}+2\irs{\Rp}{\nabla (v-u)}
\lesssim \irs{\Rp}{\nabla u}.
\end{equation}
Applying (\ref{innerregularity}) and (\ref{apriori}) to the equation (\ref{step3}), and using that $\Rp \in (\frac{R}{2}, R)$ as well as \eqref{choice} and the equality $\nabla^{tan}u=\nabla^{tan}v$ on $\partial B_{\Rp}^+$, gives that
\begin{align}
\nonumber
&\firs{r}{\nabla u - \partial_d v(0)(e_d +\nabla \tp_d)}
\\&\nonumber
\lesssim \Bigg( \left(\displaystyle\frac{r}{R}\right)^2 \firs{r}{\textrm{id}+\nabla \tp}
\\&~~~~~~~~
+\left( \displaystyle\frac{R}{r} \right)^d \left( \left(\displaystyle\frac{\rho}{R}\right)^{\beta}
\label{step4}
+\left( \displaystyle\frac{R}{\rho}\right)^{d+2} (\delta^{\h}_{R})^2 \right)  \vphantom{ \left(\displaystyle\frac{r}{R}\right)^2} \Bigg) \firs{R}{\nabla u}.
\end{align}
Now, we choose a specific $\rho$. Recall that we required  $0 < \rho \leq \frac{1}{4} \Rp$.
By varying $\rho$ subject to this condition, we can obtain $\frac{\rho}{R}=s$ for any $s  \in (0, \frac{1}{8}]$.
We select $\rho$ to satisfy $\frac{\rho}{R} = \min\{(\delta_{R}^{\h})^{2/(d+2+\beta)}, \frac{1}{8}\}$. Plugging this into (\ref{step4}) and using $\delta_R^\h\leq 1$ (which we may assume by choosing $C_\alpha(d,\lambda)$ large enough) results in 
\begin{align*}
&\firs{r}{\nabla u - \partial_d v(0)(e_d +\nabla \tp_{d})}
\\&
\lesssim \Bigg( \left(\frac{r}{R}\right)^2 \firs{r}{\textrm{id}+\nabla \tp}  
+\left( \displaystyle\frac{R}{r} \right)^d (\delta_{R}^\h)^{2\beta/(d+2+\beta)} \vphantom{\left(\frac{r}{R}\right)^2}\Bigg) \firs{R}{\nabla u}.
\end{align*}
For the first integral on the right hand side, notice that $x_d + \tp_d$ is $a$-harmonic in $B_{2r}^+$ and vanishes on $\bhs$. So, to estimate $\fint_{B_r^+} |e_d+\nabla \tp_d|^2 \,dx$ we may use (\ref{cac}).
To handle the terms of the form $e_i + \nabla \tp_i$ for $i \neq d$, we use the whole-space Caccioppoli estimate. We find that
\begin{align}
\nonumber
\firs{r}{\textrm{id}+\nabla \tp} & \lesssim \firs{r}{e_d+\nabla \tp_d} + \sum_{i=1}^{d-1} \fint_{B_r}|e_i + \nabla \tp_i|^2\,dx
\\
\label{CaccioppoliCorrector}
 & \lesssim
\frac{1}{r^2} \left( \firs{2r}{x_d + \tp_d} + \sum_{i=1}^{d-1} \fint_{B_{2r}}{|x_i+\tp_i|^2}\,dx\right).
\end{align}
Young's inequality yields
\begin{align}
\label{CaccioppoliCorrector2}
\frac{1}{r^2} \left( \firs{2r}{x_d + \tp_d}
+ \sum_{i=1}^{d-1} \fint_{B_{2r}}{|x_i+\tp_i|^2}\,dx\right) 
\lesssim \left(1 + (\delta^{\h}_{2r})^2 \right).
\end{align}
We can then conclude that
\begin{align}
\nonumber
&\firs{r}{\nabla u - \partial_d v(0) (e_d +\nabla \tp_d)}
\\&
\label{iterate1}
\lesssim  \left( \left(\frac{r}{R}\right)^2 \left(1 + \delta^2 \right)+
\left( \frac{R}{r} \right)^d \delta^{2\beta/(d+2+\beta)}\vphantom{\left(\frac{r}{R}\right)^2} \right) \firs{R}{\nabla u},
\end{align}
where we have used the notation $\delta := \max \{\delta^{\h}_{2r}, \delta^{\h}_{R}\}$.

\smallskip\smallskip
\noindent \textbf{Step 2}: Proof of the half-space excess-decay.\\
For any two radii $\tilde r$ and $\tilde R$ with $r^\ast \leq \tilde r \leq \tilde R \leq R$, we can rephrase (\ref{iterate1}) in terms of the half-space-adapted tilt-excess: Notice that for any $b \in \mathbb{R}$ the function $u - b(x_d+\tp_d)$ is $a$-harmonic on $B_{\tilde R}^+$ with homogeneous Dirichlet boundary conditions on $\bhs \cap B_R$. Applying (\ref{iterate1}) to $u- b(x_d+\tp_d)$ and taking the infimum with respect to $b$ yields
\begin{align*}
\textrm{Exc}^{\h}(\tilde r)\leq C(d, \lambda) \left( \left(\frac{\tilde r}{\overset{~}{\tilde R}}\right)^2 \left(1 + \delta^2 \right)
 + \left( \frac{\tilde R}{\tilde r} \right)^d \delta^{2\beta/(d+2+\beta)} \right) \textrm{Exc}^{\h}(\tilde R).
\end{align*}
Letting $\theta = \tilde r/\tilde R$ and using $\delta\leq 1$ gives that
\begin{equation}
 \textrm{Exc}^{\h}(\tilde r)\leq C(d, \lambda) \left(2\theta^2 + \delta^{2\beta/(d+2+\beta)} \theta^{-d} \right) \textrm{Exc}^{\h}(\tilde R),
\end{equation}
where the fixed constant $C(d, \lambda)$ comes from (\ref{iterate1}) and where we have used $\delta\leq \frac{1}{C_\alpha(d,\lambda)} \leq 1$ (the latter inequality holding w.\,l.\,o.\,g.).

We now choose $\theta$ and the constant $C_\alpha(d, \lambda)$ in the smallness condition (\ref{smallness}) in such a way that
\begin{equation}\label{k}
 C(d, \lambda)(2\theta^2+\delta^{2\beta/(d+2+\beta)}\theta^{-d})\leq \theta^{2\alpha}
\end{equation}
is satisfied. To do this we first select $\theta\in (0,1)$ such that $2C(d, \lambda) \theta^2 \leq \frac{1}{2} \theta^{2\alpha}$ holds. We then select the constant $C_\alpha(d, \lambda)$ in \eqref{smallness} to be large enough to ensure $C(d,\lambda) \delta^{2\beta/(d+2+\beta)} \theta^{-d} \leq \frac{1}{2} \theta^{2\alpha}$.
%By condition (\ref{smallness}) the choice of $C(d, \lambda, \alpha)$ then determines the minimal radius $r^*$.
This entails the estimate
\begin{align}
\label{iterate2}
\textrm{Exc}^{\h}(\theta \tilde R)\leq \theta^{2\alpha} \textrm{Exc}^{\h}(\tilde R)
\end{align}
for all $\tilde R \in [\frac{1}{\theta} r^*,R]$.

The half-space excess-decay estimate for arbitrary $r,R$ with $r^* \leq r \leq R$ follows by iterating the estimate (\ref{iterate2}). As this procedure is both straightforward and a standard argument, we omit it.

\smallskip\smallskip
\noindent \textbf{Step 3}: Proof of the coercivity of the excess expression.\\
As the left-hand side of \eqref{Coercivity} is a second-order polynomial in $b$, to establish the desired result it is sufficient to show an estimate of the form
\begin{align}
\label{coercive}
\fint_{B_r^+} |b(e_d + \nabla \tp_d)|^2 \,dx
\geq \frac{1}{2^{d+2}} |b|^2.
\end{align}
We take $\eta$ to be a cutoff with $\eta\equiv 1$ in $B_{r/2}^+$, $\eta\equiv 0$ outside $B_r^+$, $0\leq \eta\leq 1$ everywhere, and $|\nabla \eta| \leq \frac{2}{r}$. We then have
\begin{align}
\nonumber
\fint_{B_r^+} |b(e_d + \nabla \tp_d)|^2 \,dx & \geq |b|^2 \fint_{B_r^+}\eta|e_d + \nabla \tp_d|^2 \,dx
\\
\nonumber
&\geq
|b|^2
\fint_{B_r^+} \eta \,dx ~ \bigg| e_d + \frac{1}{\fint_{B_r^+} \eta \,dx} \fint_{B_r^+} \eta \nabla \phi^{\h}_{d} \,dx \bigg|^2
\\
\label{Coercive1}
&\geq
|b|^2 \fint_{B_r^+} \eta \,dx ~ \bigg| e_d - \frac{1}{\fint_{B_r^+} \eta \,dx} \fint_{B_r^+} \phi^{\h}_{d}\nabla \eta \,dx \bigg|^2.
\end{align}
Notice that the second of the above inequalities follows from an application of Jensen's inequality. Also, in the third inequality the boundary term 
has vanished due to the Dirichlet boundary conditions satisfied by $\phi^{\h}_d$.

Another use of H\"{o}lder's inequality yields that 
\begin{align*}
\displaystyle\frac{1}{\fint_{B_r^+} \eta \,dx} \bigg|\fint_{B_r^+} \phi^{\h}_{d}\nabla \eta \,dx \bigg| \leq 2^{d+1} \delta^{\h}_r. 
\end{align*}
We may assume that $C_\alpha(d,\lambda)$ in \eqref{smallness} is chosen large enough to ensure that $2^{d+1} \delta^{\h}_r \leq \frac{1}{2}$. Estimating $\fint_{B_r^+} \eta \,dx \geq ( \frac{1}{2} )^d$, we see that \eqref{coercive} now follows from \eqref{Coercive1}.

\smallskip\smallskip
\noindent \textbf{Step 4}: Proof of the mean-value property.\\
Let $r^* \leq r \leq R$; denote by $b_\rho$ the value of $b$ for which the infimum in the definition of the tilt-excess $\operatorname{Exc}^{\h}(\rho)$ is attained. We then have
\begin{align}
\nonumber
\fint_{B_r^+} |\nabla u |^2 \,dx & \lesssim \textrm{Exc}^{\h}(r) + |b_r|^2
\\& \nonumber
\lesssim \textrm{Exc}^{\h}(R) + |b_r|^2
\\&
\label{mvp2}
\lesssim \fint_{B_R^+} |\nabla u |^2 \,dx + |b_R|^2 + |b_r - b_R |^2.
\end{align}
Here, we have used \eqref{CaccioppoliCorrector}, \eqref{CaccioppoliCorrector2}, and $\delta_{2r}^\h \leq \frac{1}{C_\alpha(d,\lambda)}\leq 1$ for the first inequality, the half-space excess-decay for $\alpha = \frac{1}{2}$ for the second, and the definition of the half-space excess and Young's inequality for the third.

To complete our argument it remains to estimate $|b_R|^2$ and $|b_r-b_R|^2$. First, by the coercivity (\ref{coercive}) and the triangle inequality, we easily infer
\begin{align*}
 |b_R|^2 \lesssim \fint_{B_R^+} |b_R(e_d+\nabla \tp_d)|^2 \,dx \lesssim \textrm{Exc}^{\h}(R) + \fint_{B_R^+}| \nabla u |^2 \,dx
 \lesssim \fint_{B_R^+}| \nabla u |^2 \,dx.
\end{align*}

To estimate $|b_r-b_R|$, let $\rho\in [\max\{r^*, R/2\},R]$. Then the coercivity property (\ref{coercive}) and the triangle inequality entail
\begin{align*}
|b_{\rho}-b_{R}|^2
&\lesssim \firs{\rho}{ (b_{\rho}-b_{R})e_d + (b_{\rho}-b_{R}) \nabla \tp_d}
\\&
\lesssim \textrm{Exc}^{\h}(\rho) + \textrm{Exc}^{\h}(R)
\\&
\lesssim \fint_{B_R^+}| \nabla u |^2 \,dx.
\end{align*}

Choose $N\in \mathbb{N}_0$ such that $\frac{R}{2^{N+1}}\leq r \leq \frac{R}{2^N}$. The triangle inequality, the coercivity \eqref{coercive}, and the half-space excess-decay for $\alpha = \frac{1}{2}$ then allows us to write
\begin{align*}
|b_r - b_{R}|^2 & \leq
\left(|b_r - b_{R 2^{-N}}|+\sum_{n=1}^{N} |b_{R 2^{-n}} - b_{R 2^{-(n-1)}}|\right)^2
\\
&\lesssim \left( \sum_{n=0}^N \left( \textrm{Exc}^{\h}(R 2^{-n}) \right)^{1/2} \right)^2
\lesssim \left(\sum_{n=0}^{N} 2^{-n/2} \textrm{Exc}^{\h}(R)^{1/2}\right)^2
\\
&\lesssim \textrm{Exc}^{\h}(R).
\end{align*}
In total, (\ref{mvp2}) therefore entails the desired mean-value property.
\end{proof}

\noindent Using the half-space excess-decay we may now prove our first-order Liouville result.

\begin{proof}[Proof of Corollary \ref{Liouville}]
The Caccioppoli estimate from Lemma \ref{caclemma} shows that the growth condition (\ref{subquad}) implies that
\begin{align*}
 \lim_{R \rightarrow \infty} \displaystyle\frac{1}{R^{2\alpha}}\firs{R}{\nabla u} = 0.
\end{align*}
This, in turn, gives that
\begin{align*}
\lim_{R \rightarrow \infty} \displaystyle\frac{1}{R^{2\alpha}} \textrm{Exc}^{\h}(R)=0.
\end{align*}
By Theorem \ref{existenceofcorrectors} and Theorem \ref{excessbound} there exists a radius $r^*>0$ such that the excess-decay \eqref{ExcessDecay} holds for $R\geq r\geq r^\ast$. In particular, keeping $r$ fixed and passing to the limit $R\rightarrow\infty$, we deduce $\textrm{Exc}^{\h}(r)=0$ for any $r\geq r^\ast$. Since the coercivity property \eqref{Coercivity} implies
%we showed in Step 3 of the proof of Theorem 1
that the infimum in the definition of the excess is attained and since we have $u=0$ on $\bhs$, we find that
\begin{align*}
\forall r \geq r^* \textrm{ there exists } b\in \mathbb{R} \textrm{ such that } u(x) = b(x_d + \tp_d ) \textrm{ in } B_r^+.
\end{align*}
By the coercivity property \eqref{Coercivity}, $b$ does not depend on $r\geq r^\ast$.
Therefore, we have $u(x) = b(x_d + \tp_d)$ in $\hs$.
\end{proof}

\bibliographystyle{abbrv}
\bibliography{stochastic_homogenization}

\end{document}